%AMSTeX source file for the paper
%
%   A commuting $q$-analogue of the addition formula for 
%   disk polynomials 
%
%by P.G.A. Floris and H.T. Koelink
%Report W95-07, September 1995, Rijksuniversiteit Leiden
%%%%AMS TeX file%%%%%%%%%%%%%%%%%%%%%%%%%%%%%%%%%%%%%%%%%%%%%%%%%%%%%%%%
\magnification=1200
\input amstex
\input amssym.def
\documentstyle{amsppt}
\parindent=12pt
\baselineskip=13pt
\hsize=6.5truein
\vsize=8.9truein
%%%%%%%%%%%%%%%%%%%%%%%%%%%%Abbreviations%%%%%%%%%%%%%%%%%%%%%%%%%%%%%%%%
\def\R{{\Bbb R}}
\def\N{{\Bbb N}}
\def\C{{\Bbb C}}
\def\Z{{\Bbb Z}}

\def\Zp{{\Bbb Z}_+}
\define\z+{{\Bbb Z}_+}
\def\hf{{1\over 2}}

\def\Hi{\ell^2(\Zp)}
\def\Pol{\hbox{\it Pol}}
\def\a{\alpha}
\def\b{\beta}
\def\g{\gamma}
\def\d{\delta}
\def\e{\varepsilon}
\def\s{\sigma}
\def\t{\tau}
\def\l{\lambda}
\def\m{\mu}
\def\n{\nu}
\def\r{\rho}
\def\p{\pi}

\def\vp{\varphi}

\def\z2t{\widetilde{{\Cal Z}_2}}
\def\xt{\widetilde{\Cal X}}
\def\yt{\widetilde{\Cal Y}}

\define\ov{\overline}
\define\hb{\hfill\break}

\define\pr{{^{\prime}}}
\define\ot{\otimes}

\def\bs{\backslash}

%%%%%%%%%%%%%%%%% M a c r o s %%%%%%%%%%%%%%%%%%%%%%%%%%%%%%%%%%%%%%%%%%%
%Section, Theorem and Formula Numbering%%%%%%%%%%%%%%%%%%%%%%%%%%%%%%%%%%
%%%%%%%%%%%%%%%%%%%%%%%%%%%%%%%%%%%%%%%%%%%%%%%%%%%%%%%%%%%%%%%%%%%%%%%%%
\countdef\sectionno=1
\countdef\eqnumber=10
\countdef\theoremno=11
\countdef\countrefno=12
\countdef\cntsubsecno=13
\sectionno=0
\def\newsection{\global\advance\sectionno by 1
                \global\eqnumber=1
                \global\theoremno=1
                \global\cntsubsecno=0
                \the\sectionno}

\def\newsubsection#1{\global\advance\cntsubsecno by 1
                     \xdef#1{{\S\the\sectionno.\the\cntsubsecno}}
                     \ \the\sectionno.\the\cntsubsecno.}

\def\theoremname#1{\the\sectionno.\the\theoremno
                   \xdef#1{{\the\sectionno.\the\theoremno}}
                   \global\advance\theoremno by 1}

\def\eqname#1{\the\sectionno.\the\eqnumber
              \xdef#1{{\the\sectionno.\the\eqnumber}}
              \global\advance\eqnumber by 1}

\def\thmref#1{#1}

\global\countrefno=1

\def\refno#1{\xdef#1{{\the\countrefno}}\global\advance\countrefno by 1}

%%%%%%%%%%%%%%%%%%%%%%%%%Reference Numbering%%%%%%%%%%%%%%%%%%%%%%%%%
\refno{\AndrA}
\refno{\AskWil}
%\refno{\Brom}
\refno{\DijNo}
\refno{\Flor}
\refno{\GaspR}
\refno{\Koel}
\refno{\Koelink}
\refno{\KoelSwar}
\refno{\Koorold}
\refno{\Koor}
\refno{\Koo}
\refno{\NouMim}
\refno{\NoumYM}
\refno{\Sap}
\refno{\Stan}
\refno{\VaksS}
\refno{\AsKoo}
\refno{\Vil}
\refno{\VilKlim}

%%%%%%%%%%%%%%%%%%%%%%%%%%%%%%%%%%%%%%%%%%%%%%%%%%%%%%%%%%%%%%%%%%%%%%
%%%%%Beginning of the text%%%%%%%%%%%%%%%%%%%%%%%%%%%%%%%%%%%%%%%%%%%%
%%%%%%%%%%%%%%%%%%%%%%%%%%%%%%%%%%%%%%%%%%%%%%%%%%%%%%%%%%%%%%%%%%%%%%
\topmatter
\title A commuting $q$-analogue of the addition formula for 
disk polynomials\ \ \ \endtitle
\author P.G.A.~Floris and H.T.~Koelink\endauthor
\date Report W95-07 University of Leiden, September 1995\enddate
\address Afdeling Wiskunde en Informatica, Rijksuniversiteit Leiden,
P.O.~Box 9512, 2300~RA Leiden, the Netherlands \endaddress
\email floris\@wi.leidenuniv.nl \endemail
\address Vakgroep Wiskunde, Universiteit van Amsterdam, 
Plantage Muidergracht 24, 1018 TV Amsterdam, the Netherlands\endaddress
\email koelink\@fwi.uva.nl \endemail
\thanks Research of the second author supported by a Fellowship of the
Research Council of the Katholieke Universiteit Leuven, Belgium, and by
the Netherlands Organization for Scientific Research (NWO) under project
number 610.06.100\endthanks
\abstract  Starting from the addition formula for $q$-disk polynomials, which
is an identity in non-commuting variables, we establish a basic analogue in
commuting variables of the addition and product formula for disk
polynomials. These contain as limiting cases the addition and product formula
for little $q$-Legendre polynomials. As $q$ tends to $1$ the addition and
product formula for disk polynomials are recovered.\endabstract
\endtopmatter
\document

%\NoBlackBoxes
\noindent
1. Introduction.\hb
2. Addition formula for $q$-disk polynomials.\hb
3. Some representation theory.\hb
4. Spectral analysis of a self-adjoint operator.\hb
5. Addition formula in commuting variables.\hb
6. Product formula.\hb
7. The limit case $q\uparrow 1$.\hb
References.\hb

%%%%%%%%%%%%%%%%%%%%%%%%%%%%%%%%%%%%%%%%%%%%%%%%%%%
%%NEW SECTION %%%%%%%%%%%%%%%%%%%%%%%%%%%%%%%%%%%%%
%%%%%%%%%%%%%%%%%%%%%%%%%%%%%%%%%%%%%%%%%%%%%%%%%%%
\head\newsection . Introduction \endhead

It is known that there is a strong connection 
between the theories of
hypergeometric functions and group representations. The identification
of certain families of hypergeometric functions as specific functions on
groups enabled proofs of new identities as well
as alternative, more conceptual, proofs of well-known identities
for these families.
Particular examples of such identities are so-called addition formulas. From the 
group-theoretic point of view addition formulas reflect the homomorphism
property of group representations; for details we refer the reader to
Vilenkin \cite{\Vil} and Vilenkin and Klimyk \cite{\VilKlim}.

$q$-Special functions are generalisations (or deformations) of classical
special functions, having comparable properties, and were first studied
in the 18th and 19th century by Euler, Gauss, Heine and others. And 
although a lot is known about these functions, especially the ones of basic
hypergeometric type, there has for a long time been no satisfactory
mathematical structure on which they appear as naturally as hypergeometric
functions do on groups. The introduction of quantum groups in the 
1980's seems to have met the need for such structures.
It has now been widely accepted that quantum
groups (and quantised algebras) form a natural setting for many of
the well-known families of basic hypergeometric functions. Of particular
interest is the realisation of Askey-Wilson polynomials, at
least for specific
choices of the parameters, in quantum group theory; see
Koornwinder \cite{\Koo}, Noumi and Mimachi \cite{\NouMim},
Koelink \cite{\Koelink}, Dijkhuizen and Noumi 
\cite{\DijNo} and references given there.

The very nature of the objects one is dealing with in this area 
accounts for the appearance of identities in non-commuting
variables, or, better, in variables that satisfy certain explicit
commutation relations. Such identities are of interest in their own
right and are perhaps the natural ones to expect. However, it is
legitimate to ask whether it is possible to convert such
an identity into one involving only commuting variables without losing
any information. A way 
to do this is to represent the identity under consideration as an operator
identity in some Hilbert space and then to construct matrix elements, 
depending on certain parameters,
to obtain a scalar identity involving these additional parameters. 
It is this procedure that we use in this paper to obtain a commuting
$q$-analogue of
the addition formula for disk polynomials.

Disk polynomials are polynomials $R_{l,m}^{(\nu)}(z)$ ($l,m\in \Zp,\nu>-1$)
in the two variables $z$ and $\bar{z}$ that are orthogonal with 
respect to a positive 
measure on the closed unit disk $D$ in the complex plane. Explicitly they are
given as
$$
R_{l,m}^{(\nu)}(z) = \cases z^{l-m}
P_m^{(\nu,l-m)} (2\vert z\vert^2 - 1) &(l\geq m)\\
\bar{z}^{m-l} P_l^{(\nu,m-l)} (2\vert z\vert^2 - 1)
&(l\leq m)\endcases\tag\eqname{\defdiskpol}
$$
where the $P_k^{(\a,\b)}(x)$ denotes the Jacobi polynomial.
Their orthogonality reads
$$
\int\!\!\int_{D} R_{l,m}^{(\nu)}(re^{i\phi}) 
\ov{R_{l\pr,m\pr}^{(\nu)}(re^{i\phi})}
r (1-r^2)^{\nu}\, dr d\phi =0 \qquad \text{if } (l,m)\neq (l\pr,m\pr).
\tag\eqname{\vglorthodiskpol}
$$
For integer values $\nu=n-2$ the disk polynomials 
are the zonal spherical functions on the
compact homogeneous space $U(n)/U(n-1)$, which can be identified with
the unit sphere in $\C^n$. Here $U(n)$ denotes the group 
of unitary transformations of $\C^n$. An addition formula for these polynomials
was proved independently by ${\breve{\text S}}$apiro \cite{\Sap} and
Koornwinder \cite{\Koorold} around
1970. For $\nu>0$ it takes the following form;
$$
\eqalign{
R_{l,m}^{(\nu)}(\cos\theta_1 e^{i\phi_1} &\cos\theta_2 e^{i\phi_2}
+ \sin\theta_1 \sin\theta_2 r e^{i\psi}) = \cr
\sum_{r=0}^l\sum_{s=0}^m C_{l,m;r,s}^{(\nu)} &(\sin\theta_1)^{r+s}
 R_{l-r,m-s}^{(\nu +r+s)}(\cos\theta_1 e^{i\phi_1}) \cr 
\times &(\sin\theta_2)^{r+s}  R_{l-r,m-s}^{(\nu +r+s)}(\cos\theta_2 e^{i\phi_2})
 R_{r,s}^{(\nu -1)}(\rho e^{i\psi}),\cr}\tag\eqname{\classadd}
$$
where
$$
C_{l,m;r,s}^{(\nu)} = \frac{\nu}{\nu +r+s} {l\choose r} {m\choose s }
\frac{(\nu +l+1)_s (\nu +m+1)_r}{ (\nu +r)_s (\nu +s)_r}.
$$
The corresponding product formula, which follows
from \thetag{\classadd} and \thetag{\vglorthodiskpol}, is given by
$$
\eqalign{
&R_{l,m}^{(\nu)}(\cos\theta_1 e^{i\phi_1}) 
R_{l,m}^{(\nu)}(\cos\theta_2 e^{i\phi_2})=\cr
&\frac{\nu}{\pi} \int_0^1 \int_0^{2\pi} 
R_{l,m}^{(\nu)}(\cos\theta_1 e^{i\phi_1} \cos\theta_2 e^{i\phi_2}
+ \sin\theta_1 \sin\theta_2 r e^{i\psi})\ r (1-r^2)^{\nu -1} \, d\psi dr.\cr}
\tag\eqname{\classprod}
$$
The addition formula for 
Jacobi polynomials can be obtained from the addition formula 
\thetag{\classadd} for disk polynomials, see \cite{\Koorold} and 
\cite{\VilKlim, Vol.~2, \S 11.4.2}.

A quantum analogue of the homogeneous space $U(n)/U(n-1)$ has been 
studied by Noumi, Yamada and Mimachi \cite{\NoumYM}; see also Floris
\cite{\Flor}. They derive the following 
expression for the zonal spherical elements 
in terms of little $q$-Jacobi polynomials 
$p_k(x;q^{\a}, q^{\b};q)$ defined in (2.3); %forward reference
$$
R_{l,m}^{(\nu)}(z,z^{\ast};q^2) = \cases z^{l-m}
p_m\bigl(1-zz^{\ast};q^{2\nu},q^{2l-2m};q^2) &(l\geq m)\\
p_l\bigl(1-zz^{\ast};q^{2\nu},q^{2m-2l};q^2) (z^{\ast})^{m-l}
&(l\leq m),\endcases\tag\eqname{\defqdiskpol}
$$
where $\nu=n-2$ and $z$ and $z^{\ast}$ are the generators of a complex
$\ast$-algebra, satisfying the relation $z^{\ast}z= q^2zz^{\ast} + 1-q^2$.
Here the deformation parameter $q$ is 
assumed to be such that $0<q<1$.
The polynomials of \thetag{\defqdiskpol} are now called $q$-disk
polynomials.
 In the paper \cite{\Flor}, 
that follows the same lines as Koornwinder's treatise \cite{\Koorold}
in the classical case, 
the first author also recovers the $q$-disk polynomials as the 
zonal spherical functions on the quantum analogue of the homogeneous space
$U(n)/U(n-1)$. Moreover, an addition theorem for the $q$-disk polynomials
is proved in \cite{\Flor}, which is reproduced in section 2.
It is an identity in several
non-commuting variables which has a structure similar to \thetag{\classadd}.
It is this addition formula that we convert into one involving
only commuting variables, using the procedure mentioned before. The
outcome of it is the basic analogue of the classical addition formula for disk 
polynomials and it is the main result of this paper.\par
The remaining sections are organised as follows. 
In section 2 we 
recall some facts on $q$-disk polynomials needed later on. In 
particular we state the addition formula. Section 3
deals with irreducible $\ast$-representations of the
$\ast$-algebras on which the addition formula is realised. This information
is used in section 4 to represent the addition formula as an operator
identity on some infinite-dimensional Hilbert space. Moreover we 
give the eigenspace decomposition of this Hilbert space with respect to
a specific self-adjoint operator that appears in the left-hand side of
the addition formula under such an irreducible $\ast$-representation.  
The eigenvectors are expressible in terms of
affine $q$-Krawtchouk polynomials, and we write down explicitly the 
action of the left-hand side on these eigenvectors. In section 5 we 
calculate the image of the standard basis under the action of the
right-hand side of the operator identity. The commutative addition formula
is obtained from the fact that we use $\ast$-representations. 
In section 6 a corresponding product formula is derived.
Finally, section 7 discusses the limit transitions $q\to 1$ of both the addition
and the product formula to the classical identities 
\thetag{\classadd} and \thetag{\classprod}.

The notation for basic hypergeometric series follows the book \cite{\GaspR}
by Gasper and Rahman. By $\Zp$ we denote $\{ 0,1,2,\ldots\}$ and
$a\wedge b=\min (a,b)$ for $a,b\in\R$. Throughout this paper we fix $0<q<1$.

\demo{Acknowledgement} We thank Tom Koornwinder for his suggestions and remarks.
Part of the work of the second author was done at the Katholieke Universiteit
Leuven and thanks are due to Walter Van Assche and Alfons Van Daele for their
hospitality.
\enddemo

%%%%%%%%%%%%%%%%%%%%%%%%%%%%%%%%%%%%%%%%%%%%%%%%%%%
%%NEW SECTION %%%%%%%%%%%%%%%%%%%%%%%%%%%%%%%%%%%%%
%%%%%%%%%%%%%%%%%%%%%%%%%%%%%%%%%%%%%%%%%%%%%%%%%%%
\head\newsection . Addition formula for $q$-disk polynomials\endhead

In \cite{\Flor} an addition formula for the $q$-analogue of the disk
polynomials defined in \thetag{\defqdiskpol} has been derived. In this
section we recall this result in a slightly simpler version.
In order to present the
addition formula we first introduce certain non-commutative algebras 
and the little $q$-Jacobi polynomials that play a role in this
identity.

First we introduce the $\ast$-algebra $\xt$, which is the quotient of the
$\ast$-algebra ${\Cal X}$ of \cite{\Flor, \S 3.5} by the
ideal generated by $Q-1$. So $\xt$ is generated by $X_1$, $X_2$ subject to
the relations
$$
\gathered
X_1X_2 = qX_2X_1,\qquad X_1^\ast X_2 =qX_2X_1^\ast,\qquad
X_2^\ast X_2 = q^2 X_2X_2^\ast +1-q^2, \\
X_1^\ast X_1 = q^2 X_1X_1^\ast + (1-q^2)(1-X_2X_2^\ast).
\endgathered
\tag\eqname{\vgldefrelXtilde}
$$
We use the same notation for the generators of $\xt$ and ${\Cal X}$.

We also introduce the $\ast$-algebra $\yt$, which is the quotient of the
$\ast$-algebra ${\Cal Y}$ of \cite{\Flor, \S 3.5} by the
ideal generated by $D-1$. So $\yt$ is generated by $Y_1$, $Y_2$ subject to
the relations
$$
\gathered
Y_1Y_2 = qY_2Y_1,\qquad Y_1^\ast Y_2 =qY_2Y_1^\ast,\qquad
Y_1^\ast Y_1 = Y_1Y_1^\ast, \\
Y_1 Y_1^\ast +Y_2Y_2^\ast = 1 = q^2 Y_1^\ast Y_1 + Y_2^\ast Y_2.
\endgathered
\tag\eqname{\vgldefrelYtilde}
$$
Also here we use the same notation for the generators of $\yt$ and ${\Cal Y}$.
These $\ast$-algebras are closely related to the deformed $\ast$-algebra 
of polynomials on the sphere in $\C^n$ for $n\geq 3$ in case of $\xt$
and for $n=2$ in case of $\yt$, cf. \cite{\Flor, \S 3.5}.

The little $q$-Jacobi polynomials are defined by
$$
p_l(x;a,b;q)=\,{}_2\vp_1\left( {{q^{-l},abq^{l+1}}\atop{aq}};q,qx\right).
\tag\eqname{\vgldeflittleqJ}
$$
These polynomials are orthogonal polynomials on a geometric sequence,
see \cite{\AndrA}, \cite{\GaspR} and references therein. Using the little
$q$-Jacobi polynomials we define
$$
R^{(\n)}_{l,m}(x,y,z;q^2) = \cases z^m x^{l-m}
p_m\bigl((z-xy)z^{-1};q^{2\nu},q^{2(l-m)};q^2\bigr) &(l\geq m)\\
z^l p_l\bigl((z-xy)z^{-1};q^{2\nu},q^{2(m-l)};q^2\bigr)y^{m-l}
&(l\leq m),\endcases
\tag\eqname{\vgldefRqdisknc}
$$
which extends the definition \thetag{\defqdiskpol}. We use
\thetag{\vgldefRqdisknc} for non-commuting variables $x$, $y$, $z$, with
$yx=q^2xy +(1-q^2)z$ where
we assume that $z$ commutes with $x$ and $y$, so that \thetag{\vgldefRqdisknc} 
is polynomial in $x$, $y$ and $z$.

The following addition theorem follows directly from \cite{\Flor, Thm.~3.5.8}
with $Q=1$, $D=1$.

\proclaim{Theorem \theoremname{\thmabstractaddform}}
The following addition formula holds as an identity in $\xt\ot\yt$
for $\n>0$;
$$
\eqalign{
&R_{l,m}^{(\n)}(-qX_1\otimes Y_1^\ast+X_2\otimes Y_2
,-qX_1^\ast\otimes Y_1 + X_2^\ast\otimes Y_2^\ast ,1\otimes 1;q^2) = \cr
&\sum_{r=0}^l\sum_{s=0}^m c_{l,m;r,s}^{(\n)}(q^2)\,
R^{(\n+r+s)}_{l-r,m-s}(X_2,X_2^\ast,1;q^2)
R_{r,s}^{(\n-1)}(X_1,X_1^\ast,1-X_2X_2^\ast;q^2) \cr
&\qquad\qquad\qquad\qquad\otimes (-q)^{r-s} R^{(\n+r+s)}_{l-r,m-s}(Y_2,Y_2^\ast,1;q^2) 
Y_1^s(Y_1^\ast)^r,
\cr}
\tag\eqname{\vglabstraddform}
$$
with
$$
\align
&c^{(\n)}_{l,m;r,s}(q) = {{1-q^{\n+r+s+1}}\over{1-q^{\n+1}}}
{{c^{(\n)}_{l,m}(q)}\over{c^{(\n+r+s)}_{l-r,m-s}(q)\, c^{(\n-1)}_{r,s}(q)}},
\tag\eqname{\cees}\\
&c^{(\n)}_{l,m}(q) = {{q^{m(\n+1)}(1-q^{\n+1})}\over{1-q^{\n+l+m+1}}}
{{(q;q)_l(q;q)_m}\over{(q^{\n+1};q)_l(q^{\n+1};q)_m}}.
\endalign
$$
\endproclaim

The proof of Theorem \thmref{\thmabstractaddform} is based on the
interpretation of $q$-disk polynomials for $\n=n-2$
as zonal spherical functions on a quantum analogue of the homogeneous space
$U(n)/U(n-1)$; see also \cite{\NoumYM}. The addition formula is 
essentially the development in terms of associated spherical functions of
the outcome of the comultiplication applied to such a spherical element.

The main goal of this paper is to
deduce an addition formula in commuting
variables from \thetag{\vglabstraddform} using $\ast$-representations of
the $\ast$-algebras $\xt$ and $\yt$. These $\ast$-representations are 
studied in the next section.

%%%%%%%%%%%%%%%%%%%%%%%%%%%%%%%%%%%%%%%%%%%%%%%%%%%
%%NEW SECTION %%%%%%%%%%%%%%%%%%%%%%%%%%%%%%%%%%%%%
%%%%%%%%%%%%%%%%%%%%%%%%%%%%%%%%%%%%%%%%%%%%%%%%%%%
\head\newsection . Some representation theory \endhead

In this section we study $\ast$-representations of $\xt$ and $\yt$ using the
classification of the irreducible $\ast$-representations of the Hopf $\ast$-algebra
of polynomials on the quantum $U(n)$ group for $n=2,3$. This Hopf $\ast$-algebra
is denoted by $\Pol\,(U_q(n))$ in \cite{\Koel, \S 2}, by 
${\Cal A}_q(n)$ in
\cite{\Flor, \S 2.2} and by $A(U_q(n))$ in \cite{\NoumYM, \S 3.1}. We stick to 
the notation $\Pol\,(U_q(n))$ of \cite{\Koel}, and we also use the 
notation
$\Pol\,(SU_q(n))$ for the Hopf $\ast$-algebra which is a quotient of $\Pol\,(U_q(n))$
by the ideal generated by $D-1$. Here $D\in\Pol\,(U_q(n))$ is the 
quantum determinant \cite{\Koel, (2.7)}, which is central.

The $\ast$-algebra $\yt$ is isomorphic to $\Pol\,(SU_q(2))$ under 
the identification
$$
\pmatrix t_{11}&t_{12}\\ t_{21}&t_{22} \endpmatrix = 
\pmatrix Y_2^\ast & -qY_1^\ast\\ Y_1 & Y_2 \endpmatrix,
$$
where $t_{ij}$, $1\leq i,j\leq 2$, are the generators of $\Pol\,(SU_q(2))$
satisfying the commutation relations of \cite{\Koel, \S 2} with $D=1$.
The irreducible $\ast$-representations of this algebra have been classified
by Vaksman and Soibelman \cite{\VaksS, Thm.~3.2}; see also \cite{\Koel, Thm.~4.8}.

\proclaim{Proposition \theoremname{\propreprYt}} All mutually inequivalent
irreducible $\ast$-representations of $\yt$ are 

\parindent=8pt
\roster
\item "{{\rm{1)}}}" one-dimensional $\ast$-representations for $\phi\in[0,2\p)$;
$$
Y_1\longmapsto 0,\qquad\qquad Y_2\longmapsto e^{i\phi},
$$

\item "{{\rm{2)}}}" infinite-dimensional $\ast$-representations $\s_\phi$, 
$\phi\in[0,2\p)$, 
acting in $\Hi$ with orthonormal basis $\{ e_n\mid n\in\Zp\}$;
$$
\gather
\s_\phi(Y_1)\, e_n = e^{i\phi} q^n\, e_n, \qquad
\s_\phi(Y_2)\, e_n = \sqrt{1-q^{2n+2}}\, e_{n+1},\\
\s_\phi(Y_1^\ast)\, e_n = e^{-i\phi} q^n\, e_n, \qquad
\s_\phi(Y_2^\ast)\, e_n = \sqrt{1-q^{2n}}\, e_{n-1},
\endgather
$$
where $e_n=0$ for negative $n$ by convention.
\endroster
\parindent=12pt
\endproclaim

\demo{Remark \theoremname{\remtoproponirrrepY}} Although $\s_\phi$ is not
faithful, it is known that $\cap_\phi {\text{ker}}\, (\s_\phi)$ is trivial.
\enddemo

To describe irreducible $\ast$-representations of $\xt$ suited for our purposes we 
use the embedding $\xt\hookrightarrow\Pol\,(U_q(3))$ as $\ast$-algebras given by
$$
X_1\longmapsto t_{21}^\ast,\qquad X_2\longmapsto t_{11}^\ast.
\tag\eqname{\vglembeddingXt}
$$
The commutation relations \thetag{\vgldefrelXtilde} for $\xt$ correspond
to \cite{\Koel, (2.1), (2.14), (2.16)}. Note that \thetag{\vglembeddingXt}
also gives an embedding $\xt\hookrightarrow\Pol\,(U_q(n))$ for $n\geq 3$
(this embedding also follows from a combination of \cite{\Flor, Prop.~3.5.7}
with \cite{\Flor, Prop.~3.1.4} and the fact that the map $t_{ij}\to
t_{n-j+1,n-i+1}$ extends to an algebra anti-automorphism of
$\Pol\, (U_q(n))$).
The irreducible $\ast$-representations of $\Pol\,(U_q(n))$ have been classified
in \cite{\Koel, Thm.~4.10}, and by specialising $n=3$ and restricting
to $\xt$ we obtain families of $\ast$-representations of $\xt$. Note that 
we do not claim the list is exhaustive.

\proclaim{Proposition \theoremname{\propreprXt}} The following list gives 
mutually inequivalent irreducible $\ast$-re\-pre\-sen\-ta\-tions of $\xt$;

\parindent=8pt
\roster
\item "{{\rm{1)}}}" one-dimensional $\ast$-representations for $\phi\in[0,2\p)$
$$
X_1\longmapsto 0, \qquad\qquad X_2\longmapsto e^{i\phi},
$$

\item "{{\rm{2)}}}" infinite-dimensional $\ast$-representations $\p_1^\phi$, 
$\phi\in[0,2\p)$, acting on
$\Hi$ with orthonormal basis $\{ e_n \mid n\in\Zp\}$;
$$
\p_1^\phi(X_1)\, e_n = e^{i\phi} q^n\, e_n, \qquad
\p_1^\phi(X_2)\, e_n = \sqrt{1-q^{2n+2}}\, e_{n+1},
$$

\item "{{\rm{3)}}}" infinite-dimensional $\ast$-representation $\p$ acting on 
$\ell^2(\Zp^2)$ with orthonormal basis $\{f_\m \mid \m=(\m_1,\m_2)\in\Zp^2\}$;
$$
\gather
\p(X_1)\, f_\m = q^{\m_2} \sqrt{1-q^{2\m_1+2}}\, f_{\m+\e_1},\qquad
\p(X_2)\, f_\m = \sqrt{1-q^{2\m_2+2}}\, f_{\m+\e_2},\\
\p(X_1^\ast)\, f_\m = q^{\m_2} \sqrt{1-q^{2\m_1}}\, f_{\m-\e_1},\qquad
\p(X_2^\ast)\, f_\m = \sqrt{1-q^{2\m_2}}\, f_{\m-\e_2},
\endgather
$$
with $\e_1=(1,0)$, $\e_2=(0,1)$ and the convention that $f_\m=0$ if $\m_1$ or $\m_2$
is negative. This representation is a faithful representation of $\xt$.
\endroster
\parindent=12pt
\endproclaim

\demo{Proof} Notice that the algebra $\yt$ is isomorphic to the algebra
$\xt$ divided out by the ideal generated by $X_1^{\ast}X_1 - X_1 X_1^{\ast}$.
So parts 1) and 2) follow from Proposition \thmref{\propreprYt}.
It is obvious that the mapping $\pi$ in 3) defines an irreducible
$\ast$-representation of $\xt$ which is not equivalent to the ones in
1) and 2). Hence it remains to prove that $\pi$ is faithful.

Let $\xi\in\xt$ satisfy $\p(\xi)=0$ and write $\xi$ as a finite sum;
$$
\xi=\sum_{r,s,t,u\in\Zp} c_{r,s,t,u}\, X_1^r(X_1^\ast)^s X_2^t(X_2^\ast)^u.
$$
To simplify calculations we use the orthogonal basis
$g_\m= \left( (q^2;q^2)_{\m_1}(q^2;q^2)_{\m_2}\right) ^{\hf}\, f_\m.$
Then for all $\m\in\Zp^2$,
$$
0 = \p(\xi)\, g_\m = 
\sum_{r,s,t,u\in\Zp} c_{r,s,t,u}\, (q^{2\m_1};q^{-2})_s (q^{2\m_2};q^{-2})_u
q^{(\m_2+t-u)(s+r)} \, g_{\m+(r-s)\e_1 + (t-u)\e_2}.
$$
Consider the coefficient of $g_{\m+\a\e_1+\b\e_2}$ to find, for fixed $\a$, $\b$,
$$
\sum_{r=s+\a}\sum_{t=u+\b}  c_{r,s,t,u}\, (q^{2\m_1};q^{-2})_s (q^{2\m_2};q^{-2})_u
q^{(\m_2+t-u)(s+r)} = 0, \quad\forall\, \m_1,\m_2\in\Zp.
$$
This double sum is a polynomial in two variables, $q^{\m_1}$, $q^{\m_2}$, 
so we find that $c_{r,s,t,u}=0$ for $r=s+\a$, $t=u+\b$. 
Since $\a$, $\b$ are arbitrary, 
all coefficients are zero. Thus $\xi=0$.
\qed\enddemo

As a consequence of the proof of Proposition \thmref{\propreprXt} we
obtain

\proclaim{Corollary \theoremname{\corpropreprX}}  
The set $\{ X_1^r(X_1^\ast)^s X_2^t(X_2^\ast)^u\mid r,s,t,u\in\Zp\}$ forms 
a linear basis for $\xt$.
\endproclaim

\demo{Remark \theoremname{\remtopropreprX}} 
Instead of checking the commutation relations in the proof of
Proposition \thmref{\propreprXt} we can also
observe that the actions are coming from the restrictions of certain
irreducible $\ast$-representations of $\Pol\,(U_q(3))$. According to
\cite{\Koel, Thm.~4.10} the irreducible $\ast$-representations of $\Pol\,(U_q(3))$
are parametrised by $\rho\in S_3$ and $(\g_1,\g_2,\g_3)\in\C^3$ with
$|\g_i|$ determined by $\rho$. In case 1) the $\ast$-representation follows from
restriction of the irreducible $\ast$-representation of $\Pol\,(U_q(3))$
corresponding to $\rho=1$, $(\g_1,\g_2,\g_3)=(e^{-i\phi},1,1)$; in case 2)
we take $\rho=(12)$, $(\g_1,\g_2,\g_3)=(e^{-i\phi},-q,1)$ and in case 3)
we take $\rho=(123)$,  $(\g_1,\g_2,\g_3)=(1,-q,-q)$. The orthonormal bases
in the cases 2) and 3) correspond to the orthogonal bases described
in \cite{\Koel, Cor.~4.15(i)}.
\enddemo

%%%%%%%%%%%%%%%%%%%%%%%%%%%%%%%%%%%%%%%%%%%%%%%%%%%
%%NEW SECTION %%%%%%%%%%%%%%%%%%%%%%%%%%%%%%%%%%%%%
%%%%%%%%%%%%%%%%%%%%%%%%%%%%%%%%%%%%%%%%%%%%%%%%%%%
\head\newsection . Spectral analysis of a self-adjoint operator \endhead

Application of a one-dimensional $\ast$-representation of either $\xt$ or
$\yt$, which are described in \S 3, 
to the addition formula \thetag{\vglabstraddform} leads to a trivial
identity, so in order to convert the addition formula \thetag{\vglabstraddform} 
into an addition formula in commuting variables we apply the infinite-dimensional
$\p\otimes \s$ to the identity \thetag{\vglabstraddform} in
$\xt\otimes\yt$, where $\p$ and $\s=\s_0$ are irreducible $\ast$-representations
of $\xt$ and $\yt$ introduced in the previous section.
In this section we study the operator emerging from the
left-hand side of \thetag{\vglabstraddform}.

Consider the elements $T, T^{\ast}\in\xt\otimes\yt$ defined by
$$
T = -q X_1\otimes Y_1^\ast + X_2\otimes Y_2,\qquad
T^{\ast} = -q X_1^{\ast}\otimes Y_1 + X_2^{\ast}\otimes Y_2^{\ast}
\tag\eqname{\vgldefT}
$$
and the positive self-adjoint element $R=TT^\ast$. Then the left-hand side
of \thetag{\vglabstraddform} is a polynomial in $1-R$ multiplied from the
left, respectively right, by some power of $T$, respectively $T^\ast$.
Here $1=1\otimes 1$ is the identity of $\xt\otimes\yt$.

Using Propositions \thmref{\propreprYt} and {\propreprXt} we first calculate
$$
\eqalign{
\bigl( (\p\otimes\s)T\bigr) f_\m\otimes e_k =
\sqrt{(1-q^{2\m_2+2})(1-q^{2k+2})} f_{\m+\e_2}&\otimes e_{k+1} \cr
 - q^{1+\m_2+k} \sqrt{1-q^{2\m_1+2}} &f_{\m+\e_1}\otimes e_k\cr}
\tag\eqname{\vglactTbv}
$$
and
$$
\eqalign{
\bigl( (\p\otimes\s)T^\ast\bigr) f_\m\otimes e_k =
\sqrt{(1-q^{2\m_2})(1-q^{2k})} f_{\m-\e_2}&\otimes e_{k-1}\cr
 - q^{1+\m_2+k} \sqrt{1-q^{2\m_1}} &f_{\m-\e_1}\otimes e_k.\cr}
\tag\eqname{\vglactTstarbv}
$$

By $H(\m,p)$ we denote the $N+1$-dimensional subspace
of the representation space $\ell^2(\Zp^2)\otimes\Hi$ spanned by the
vectors $f_{\m+(l-p\wedge\m_2)\e}\otimes e_{p+l-p\wedge\m_2}$,
$l=0,\ldots,N$ with $N=\m_1+p\wedge\m_2$ and $\e=\e_2-\e_1=(-1,1)$.
Note that all these finite-dimensional spaces are orthogonal;
$H(\m,p)\perp H(\l,m)$ if $(\m,p)\not= (\l,m)$ and
$\ell^2(\Zp^2)\otimes\Hi=\bigoplus_{\m,p} H(\m,p)$.
It follows from \thetag{\vglactTbv} and \thetag{\vglactTstarbv} that
$$
(\p\otimes\s)T\ \colon\ H(\m,p)\to H(\m+\e_1,p),\quad
(\p\otimes\s)T^\ast\ \colon\ H(\m,p)\to H(\m-\e_1,p).
$$

Then $(\p\otimes\s)R\ \colon\ H(\m,p)\to H(\m,p)$ and we show
how to diagonalise this finite-dimensional mapping.
From \thetag{\vglactTbv} and \thetag{\vglactTstarbv} we obtain,
with $\e=\e_2-\e_1$ as before,
$$
\eqalign{
\bigl( (\p\otimes\s)R\bigr)(f_\m\otimes e_k) =
&\bigl( (1-q^{2\m_2})(1-q^{2k})
+q^{2+2k+2\m_2}(1-q^{2\m_1})\bigr)(f_\m\otimes e_k)\cr
-&q^{k+\m_2+1}\sqrt{(1-q^{2\m_1})(1-q^{2\m_2+2})(1-q^{2k+2})}
(f_{\m+\e}\otimes e_{k+1})\cr
-&q^{k+\m_2-1}\sqrt{(1-q^{2\m_2})(1-q^{2\m_1+2})(1-q^{2k})}
(f_{\m-\e}\otimes e_{k-1}).\cr}\tag\eqname{\vglactionR}
$$
Hence $(\p\otimes\s)R$ on the finite-dimensional space $H(\m,p)$
is given by a Jacobi matrix, i.e. a tridiagonal symmetric matrix, so
the eigenvectors are described in terms of discrete orthonormal
polynomials on a finite set. Let us first
assume that $p\geq \m_2$. It follows from \thetag{\vglactionR} that
the vector $\sum_{j=0}^{\m_1+\m_2} p_j(\l) f_{\m+(j-\m_2)\e}\otimes
e_{p+j-\m_2}$
is an eigenvector for $(\p\otimes\s)R$ for the eigenvalue $\l$ if and
only if
$$
\eqalign{
&(\l+q^{2N}-1) p_j(\l) = a_j p_{j+1}(\l) + b_j p_j(\l) + a_{j-1} p_{j-1}(\l),\cr
&(\l+q^{2N}-1) p_N(\l) = b_N p_N(\l) + a_{N-1} p_{N-1}(\l)\cr}
$$
with
$$
\align
a_j &= - q^{2j}\sqrt{q^{2+2y}(1-q^{2j+2})(1-q^{2N-2j})(1-q^{2j+2y+2})}, \\
b_j &= q^{2N}\bigl( (1-q^{2j-2N})(1-q^{2y+2+2j}) -
q^{2y+2j-2N}(1-q^{2j})\bigr),
\endalign
$$
where $y=p-\m_2$, $N=\m_1+\m_2$ and the standard initial conditions
$p_{-1}(\l)=0$, $p_0(\l)=1$ hold. Then $p_j(\l)$
is a polynomial of degree $j$ in $\l$, $0\leq j\leq N$.

The polynomial $p_j(\l)$ can be written in terms of the affine
$q$-Krawtchouk polynomial as follows. Recall that the
affine $q$-Krawtchouk polynomial is defined by
$$
K_l(x) =\, 
K_l(x;t,N;q) =\, {}_3\vp_2\left( {{q^{-l},x,0}\atop{tq,q^{-N}}};q,q\right),
\qquad 0\leq l\leq N,
\tag\eqname{\vgldefaffqK}
$$
and satisfies the three-term recurrence relation
$$
\multline
(1-q^{-x})K_l(q^{-x}) = (q^{l-N}-1)(1-tq^{l+1})K_{l+1}(q^{-x}) \\
+\Bigl( (1-q^{l-N})(1-tq^{n+1})-tq^{l-N}(1-q^l)\Bigr) K_l(q^{-x}) 
+tq^{l-N}(1-q^l)K_{l-1}(q^{-x}),
\endmultline
$$
cf. Stanton \cite{\Stan, \S 4}, Askey and Wilson \cite{\AskWil} and references
therein. Denote by
$$
\hat K_l (x) = \hat K_l(x;t,N;q) =
(-1)^l (tq)^{-l/2} \biggl( {{(q^N;q^{-1})_l
(tq;q)_l} \over{(q;q)_l}} \biggr)^\hf K_l(x;t,N;q)
\tag\eqname{\vgldefortonormalaffqK}
$$
the orthonormal affine $q$-Krawtchouk polynomials with positive
leading coefficient, where we assume $0<tq<1$. Then
$$
\eqalign{
&(1-q^{-x})\hat K_l(q^{-x}) = a_l \hat K_{l+1}(q^{-x}) +
b_l \hat K_l(q^{-x}) + a_{l-1} \hat K_{l-1}(q^{-x}),\cr
&(1-q^{-x})\hat K_N(q^{-x}) = 
b_N \hat K_N(q^{-x}) + a_{N-1} \hat K_{N-1}(q^{-x})\cr}
\tag\eqname{\recurrorthonormKraw}
$$
for $x\in\{ 0,\ldots, N\}$ with
$$
\align
a_l &= a_l(t,N;q)= - q^{l-N}\sqrt{tq (1-tq^{l+1}) 
(1-q^{l+1}) (1-q^{N-l})},\\
b_l &= b_l(t,N;q)= (1-q^{l-N}) (1- tq^{l+1}) - tq^{l-N}(1-q^l).
\tag\eqname{\recurrcoeff}\endalign
$$
Comparing both three-term recurrence relations shows that we can make the choice
$p_j(\l)= \hat K_j(q^{-2x};q^{2y},N;q^2)$ with
$\l=1-q^{2(N-x)}$ for $x\in\Zp$, $0\leq x\leq N$ and
$y=p-\m_2$, $N=\m_1+\m_2$ as before.

In case $p<\m_2$ we look for an
eigenvector of the form  $\sum_{j=0}^{p+\m_1} p_j(\l) f_{\m+(j-p)\e}\otimes
e_{j}$ and we obtain the
same recurrence relation for the polynomials $p_j(\l)$ as above,
but with $y=\m_2-p$ and $N=\m_1+p$. This follows from the symmetry of 
\thetag{\vglactionR} in $k$ and $\m_2$.

\demo{Remark \theoremname{\remorthoaffqKrawtchouk}}
We have $K_j(q^{-x};t,N;q)=K_x(q^{-j};t,N;q)$, $0\leq j,x\leq N$, which means
that the affine $q$-Krawtchouk polynomials are self-dual. Using the
self-duality and the fact that they form eigenvectors of the self-adjoint
operator $(\p\otimes\s)R$ (with other parameters) for different eigenvalues 
gives the
known orthogonality relations \cite{\AskWil}, \cite{\Stan}
$$
\sum_{j=0}^N {{(q^N;q^{-1})_j(tq;q)_j}\over{(q;q)_j}} (tq)^{-j} K_x(q^{-j})
K_y(q^{-j}) =
\d_{xy} {{(q;q)_x}\over{(q^N;q^{-1})_x (tq;q)_x}}
(tq)^{x-N}
\tag\eqname{\Kraworth}
$$
for the affine $q$-Krawtchouk polynomials
$K_x(q^{-j})=K_j(q^{-x})=K_j(q^{-x};t,N;q)$. This follows from the previous
discussion apart for the occurrence of $(tq)^{-N}$ on the right-hand side of
\thetag{\Kraworth}. It can be evaluated by taking $x=y=0$ and
reversing the summation parameter. The resulting sum can be calculated using
the $q$-Chu-Vandermonde sum. The
corresponding orthogonality measure is positive definite if $0<tq<1$.
\enddemo

\proclaim{Proposition \theoremname{\propeigvectReen}}
For the $N+1=\m_1+p\wedge\m_2 +1$-dimensional subspace $H(\m,p)$,
$\m\in\Zp^2$, $p\in\Zp$, of 
the representation space $\ell^2(\Zp^2)\otimes\Hi$ 
corresponding to the representation
$\p\otimes\s$ of $\xt\otimes\yt$ there exists an orthogonal basis 
$\{ g_x\mid x\in\{ 0,\ldots, N\}\}$ of eigenvectors for the positive operator
$(\p\otimes\s)R$ (with $R=TT^\ast \in\xt\otimes\yt$
defined by \thetag{\vgldefT}). Explicitly, 
$\bigl( (\p\otimes\s)R\bigr)g_x= (1-q^{2x}) g_x$ and
$$
g_x = g_x(\m,p) = \sum_{j=0}^{N} \hat K_j(q^{2x-2N};q^{2y},N;q^2)
f_{\m+(j-p\wedge\m_2)\e}\otimes e_{p+j-p\wedge\m_2}
$$
with $y=|p-\m_2|$ and $\e=\e_2-\e_1$. The orthogonality relations are
$$
\langle g_x, g_{x\pr}\rangle = \d_{xx\pr}
{{(q^2;q^2)_{N-x} q^{-2x(y+1)}}\over{(q^{2N};q^{-2})_{N-x}
(q^{2y+2};q^2)_{N-x} }} .
$$
\endproclaim

\demo{Remark \theoremname{\remthmpropeigvectReen}}  
Instead of taking the representation $\p$ of $\xt$ we can take the 
representation $\p_1^\phi$ of $\xt$ described in Proposition 
\thmref{\propreprXt}. Then we can find in a similar way eigenvectors
of $(\p_1^\phi\otimes\sigma)R$ in terms of Wall polynomials.
Define the orthonormal Wall polynomials by
$$
w_l(x;t;q) = (-1)^l (tq)^{-l/2} \biggl({{(tq;q)_l}\over{(q;q)_l}}\biggr)^\hf
\, {}_2\vp_1\left( {{q^{-l},0}\atop{tq}};q,qx\right);
\tag\eqname{\vgldefWallorthon}
$$
see e.g. \cite{\Koor}, \cite{\AsKoo} and references therein.
Note that the Wall polynomials are little $q$-Jacobi polynomials 
\thetag{\vgldeflittleqJ} with $b=0$.
The eigenvectors $g_x$ of $(\p_\phi^1\otimes\s)R$ for the
eigenvalue $1-q^{2x}$, $x\in\Zp$, in $\Hi\otimes\Hi$ are
$$
g_x = \sum_{j=0}^\infty e^{-ij\phi} w_j(q^{2x};q^{2p};q^2) 
e_j\otimes e_{p+j},\quad
g_x = \sum_{j=0}^\infty e^{-ij\phi} w_j(q^{2x};q^{2p};q^2) e_{p+j}\otimes e_j,
\tag\eqname{\vgldefeigvectWall}
$$
for $p\in\Zp$. For $\phi=0$ these eigenvectors can be obtained by a formal
limiting process of the result of Proposition \thmref{\propeigvectReen}.
Take $\m_2=0$ and let $\m_1\to\infty$ and use the limit transition
$$
\lim_{N\to\infty} \hat K_l(xq^{-N};t,N;q) = w_l(x;t;q)
\tag\eqname{\vgllimaffqKtoWall}
$$ 
of the affine $q$-Krawtchouk polynomials to the Wall polynomials. 
So here we let $e_k\in\Hi$ correspond to $\lim_{\m_1\to\infty} f_{\m_1,k}$
in $\ell^2(\Zp^2)$. It follows from Proposition \thmref{\propreprXt} that
under this limit the $\ast$-representation $\p$ formally goes over into $\pi_1^0$.

The occurrence of the Wall polynomials in relation to the eigenvectors of
$(\p_0^1\otimes\s)R$ as in \thetag{\vgldefeigvectWall} is contained
in Koornwinder \cite{\Koor}. To see this, recall, cf. proof of 
Proposition \thmref{\propreprXt}, that $\p^1_0$ is a representation
of $\xt/ <X_1^\ast X_1 - X_1^\ast X_1> \cong \yt$, and under this isomorphism
$\p^1_0=\s$. Now $\yt$ is isomorphic to the $\ast$-algebra ${\Cal A}$ of
\cite{\Koor, \S 3} by $Y_1\mapsto\g$, $Y_2\mapsto \a^\ast$ and $\p^1_0=\s$
then corresponds to the representation $\t$ of \cite{\Koor, (3.6)}.
Under this identification $T$ corresponds to $\Phi(\a^\ast)$ of 
\cite{\Koor, \S 3}, where $\Phi$ is the notation for the
comultiplication in \cite{\Koor}. Hence $1-R$ corresponds to
$1-\Phi(\a^\ast)\Phi(\a)=\Phi(1-\a^\ast\a)=\Phi(\g\g^\ast)$, and
now \cite{\Koor, (4.6)} corresponds precisely to the vectors in
\thetag{\vgldefeigvectWall} being eigenvectors of $(\p^1_\phi\otimes\s)R$
for $\phi=0$.

The orthogonality relations for $g_x$, i.e. the orthogonality relations
dual to the orthogonality relations for the Wall polynomials, now correspond
to the orthogonality relations for the Al-Salam--Carlitz polynomials
$V_n^{(a)}$. Explicitly,
$$
%\gather
w_l(q^n;t;q) = {{t^{l/2}q^{l^2/2}}\over{\sqrt{(tq,q;q)_l}}} 
\, {}_2\vp_0\left( {{q^{-l},q^{-n}}\atop{-}};q,{{q^n}\over t}\right) %\\
= {{(-t)^{-n}q^{n(n-1)/2}t^{l/2}q^{l^2/2}}\over{\sqrt{(tq,q;q)_l}}}
V_n^{(t)}(q^{-l};q),
%\endgather
$$
where the first equality follows from \cite{\GaspR, (III.6), (III.7)}.
This is implicitly contained in \cite{\Koor, \S 2}, but the connection
with Al-Salam--Carlitz polynomials is not noticed.
\enddemo

To investigate how the mappings $(\p\otimes\s)T$ and
$(\p\otimes\s)T^\ast$ act on the eigenvectors described in
Proposition \thmref{\propeigvectReen}, we first prove a lemma.

\proclaim{Lemma \theoremname{\lemcomrelRT}}
In $\xt\otimes\yt$ we have the commutation relations
$$
(1-R)T=q^2T(1-R),\qquad T^\ast (1-R) = q^2 (1-R) T^\ast.
$$
\endproclaim

\demo{Proof} Since $R=TT^\ast$ is self-adjoint, it suffices to prove one of the
commutation relations. The relation $(1-R)T=q^2T(1-R)$ is implied
by $T^\ast T-q^2TT^\ast=1-q^2$, which is proved directly from 
\thetag{\vgldefT} and the
commutation relations \thetag{\vgldefrelXtilde} and \thetag{\vgldefrelYtilde}
for $\xt$ and $\yt$.
\qed\enddemo

\demo{Remark \theoremname{\remlemcomrelRT}}
As stated in section 2 the algebras $\xt$ and $\yt$ are closely related to
a deformation of the $\ast$-algebra of polynomials on the sphere in $\C^n$
for $n=3$ and $n=2$. This deformation is denoted by $\widetilde{{\Cal Z}_n}$
in \cite{\Flor} on which $\Pol\,(U_q(n))$ acts in a natural way. It is
possible to show that $(1-R)T = q^2 T (1-R)$ follows from the commutation 
relation $Q_{n-1}z_n=q^2z_nQ_{n-1}$ in $\widetilde{{\Cal Z}_n}$ in the notation
of \cite{\Flor} under a suitable algebra homomorphism. 
This leads to a more conceptual proof of Lemma \thmref{\lemcomrelRT}, but
falls outside the scope of this paper.
\enddemo

Lemma \thmref{\lemcomrelRT} implies that
$\bigl((\p\otimes\s)T \bigr)g_x(\m,p)$
is an eigenvector of $(\p\otimes\s)(1-R)$ for the eigenvalue
$q^{2x+2}$ in the space $H(\m+\e_1,p)$. Since the eigenvectors
$g_x(\m+\e_1,p)$, $x=0,\ldots,\m_1+p\wedge\m_2+1$, form an orthogonal basis of
$H(\m+\e_1,p)$, we have proved part of the following proposition.

\proclaim{Proposition \theoremname{\propactieTev}}
For $\m\in\Zp^2$, $p\in\Zp$ we have
$$
\align
\bigl( (\p\otimes\s)T \bigr)g_x(\m,p) &= -q^{1+y}
\sqrt{1-q^{2N+2}}\, g_{x+1}(\mu+\e_1,p), \\
\bigl( (\p\otimes\s)T^\ast \bigr) g_x(\m,p) &= {{-q^{-y-1}}\over
{\sqrt{1-q^{2N}}}} (1-q^{2x})\, g_{x-1}(\mu-\e_1,p),
\endalign
$$
where $T,T^{\ast}\in\xt\otimes\yt$ are defined by \thetag{\vgldefT}, 
$y=\vert p-\m_2\vert$, $N=\m_1+p\wedge\m_2$
and the vectors $g_x(\m,p)\in\ell^2(\Zp^2)\otimes\Hi$ are defined in Proposition
\thmref{\propeigvectReen}.
\endproclaim

\demo{Proof} We have that $\bigl((\p\otimes\s)T\bigr) g_x(\m,p)=C
g_{x+1}(\mu+\e_1,p)$ for some constant $C$, as already remarked.
To calculate $C$, take the inner product with
$f_{(N+1,\m_2- p\wedge\m_2)}\otimes e_{p-p\wedge\m_2}
\in H(\m+\e_1,p)$
and use the fact that we are dealing with a $\ast$-representation to get
$$
\align
C &= \langle C g_{x+1}(\m+\e_1,p), 
f_{(N+1,\m_2-p\wedge\m_2)}\otimes e_{p-p\wedge\m_2}
\rangle \\
&= \langle g_x(\m,p), \bigl( (\p\otimes\s)T^\ast\bigr)
f_{(N+1,\m_2-p\wedge\m_2)}\otimes e_{p-p\wedge\m_2}
\rangle \\
&= -q^{1+y} \sqrt{1-q^{2N+2}},
\endalign
$$
where we use \thetag{\vglactTstarbv} to obtain the last equality.
We also have for some other constant $C$ that
$\bigl( (\p\otimes\s)T^\ast\bigr) g_x(\m,p)=C g_{x-1}(\mu-\e_1,p)$.
Now use this and the first statement of the proposition to get
$$
q^{2x} g_x(\m,p) = \bigl( (\p\otimes\s)(1-T T^\ast)\bigr) g_x(\m,p)
= (1+C q^{1+y}\sqrt{1-q^{2N}})\, g_x(\m,p)
$$
from which $C$ follows.
\qed\enddemo

\demo{Remark \theoremname{\remcontigrelaffqK}} 
From Proposition \thmref{\propactieTev}
recurrence formulas for affine
$q$-Krawtchouk polynomials with shifted parameters can be obtained.
Combining Proposition \thmref{\propactieTev} with
\thetag{\vglactTbv} and replacing $q^2$, $x-N$, $q^{2y}$ by $q$, $-x$, $t$
we get
$$
\eqalign{
(1-q^{N+1}) K_l(q^{-x}; t, N+1;q) = 
q^l (1-&q^{N-l+1}) K_l(q^{-x}; t, N;q) \cr
&+(1-q^l) K_{l-1}(q^{-x}; t, N;q).
\cr}
$$
Similarly, combination of
Proposition \thmref{\propactieTev} with \thetag{\vglactTstarbv} gives
$$
\eqalign{
(1-q^{N-x}) K_l(q^{-x}; t, N-1;q) = 
tq^{l+1} &(1-q^N) K_l(q^{-x}; t, N;q)\cr
&+(1-q^N)(1-tq^{l+1}) K_{l+1}(q^{-x}; t, N;q).
\cr}
$$
These relations are contiguous relations for the ${}_3\vp_2$-series in
\thetag{\vgldefaffqK}.
\enddemo

The next corollary shows how the adjoint of the left-hand side of the addition
formula \thetag{\vglabstraddform} under $\p\otimes\s$ acts on the eigenvectors
of Proposition \thmref{\propeigvectReen}.
In order to formulate the corollary we introduce the following notation;
$$
r^{(\n)}_{l,m}(x;q) =\cases (xq;q)_{l-m}^\hf p_m(x;q^\n,q^{l-m};q)
&(l\geq m)\\
(x;q^{-1})_{m-l}^\hf p_l(xq^{l-m};q^\n,q^{m-l};q) &(l\leq m),
\endcases
\tag\eqname{\vglqdiskpolwithoutphase}
$$
in terms of the little $q$-Jacobi polynomial \thetag{\vgldeflittleqJ}.
Note that $r^{(\n)}_{l,m}(xq^{m-l};q) = r^{(\n)}_{m,l}(x;q)$.

\proclaim{Corollary \theoremname{\coractiolhs}}
With the notation as in \thetag{\vgldefRqdisknc}, \thetag{\vgldefT},
Proposition \thmref{\propeigvectReen} and 
with $N=\m_1+p\wedge\m_2,\ y=\vert p-\m_2\vert$ we have
$$
\eqalign{
\Bigl( (\p\otimes\s)R^{(\n)}_{l,m}(T,T^\ast,1;q^2) \Bigr)^\ast\ g_x(\m,p) =\ 
& (-q^{y+1})^{m-l} (q^{2N+2};q^2)_{m-l}^\hf 
\ (q^{2x-2l+2m+2};q^2)_{l-m}^\hf\cr
&\times r_{l,m}^{(\n)}(q^{2x-2l+2m};q^2)\ g_{x+m-l}(\m+(m-l)\e_1,p).
\cr}
$$
Here we use the standard notation \cite{\GaspR} 
$(a;q)_{n}=(a;q)_\infty/(aq^{n};q)_\infty$ for $n\in\Z$.
\endproclaim

%%%%%%%%%%%%%%%%%%%%%%%%%%%%%%%%%%%%%%%%%%%%%%%%%%%
%%NEW SECTION %%%%%%%%%%%%%%%%%%%%%%%%%%%%%%%%%%%%%
%%%%%%%%%%%%%%%%%%%%%%%%%%%%%%%%%%%%%%%%%%%%%%%%%%%
\head\newsection . Addition formula in commuting variables \endhead

In the previous section we have studied the left-hand side of the 
addition formula \thetag{\vglabstraddform} under the $\ast$-representation
$\p\otimes\s$ of $\xt\otimes\yt$. In this section we first study the right
hand side of the 
addition formula \thetag{\vglabstraddform} under the $\ast$-representation
$\p\otimes\s$ of $\xt\otimes\yt$. This then directly leads to an addition
formula for the functions $r_{l,m}^{(\n)}(\cdot;q)$ as defined in
\thetag{\vglqdiskpolwithoutphase}.

Firstly, from \thetag{\vgldefRqdisknc} and Proposition \thmref{\propreprYt} 
we see
that the action of $\s$ on the second leg of the tensor product
of the right-hand side of the  addition formula
\thetag{\vglabstraddform} is given by
$$
\s\bigl(R^{(\n+r+s)}_{l-r,m-s}(Y_2,Y_2^\ast,1;q^2) Y_1^s (Y_1^\ast)^r\bigr) e_k
= q^{k(r+s)} r^{(\n+r+s)}_{l-r,m-s}(q^{2k};q^2)\ e_{k+l-m+s-r}.
\tag\eqname{\vglactionsltp}
$$
Next we calculate the action of $\p$ on the first leg of the
tensor product of the right
hand side of the addition formula \thetag{\vglabstraddform}.
By iteration of Proposition \thmref{\propreprXt} we see that
$$
\gather
\p(X_1^p)f_\m = q^{p\m_2} (q^{2\m_1+2};q^2)_p^\hf f_{\m+p\e_1},\quad
\p(X_2^p)f_\m = (q^{2\m_2+2};q^2)_p^\hf f_{\m+p\e_2},\\
\p\bigl((X_1^\ast)^p\bigr)f_\m = q^{p\m_2} (q^{2\m_1};q^{-2})_p^\hf
f_{\m-p\e_1},\quad
\p\bigl((X_2^\ast)^p\bigr)f_\m = (q^{2\m_2};q^{-2})_p^\hf f_{\m-p\e_2}.
\endgather
$$
Consequently, using $\p(1-X_2X_2^\ast)f_\m = q^{2\m_2}f_\m$ we get
$$
\p\bigl(R_{r,s}^{(\n-1)}(X_1,X_1^\ast,1-X_2X_2^\ast;q^2)\bigr) f_\m =
 q^{\m_2(r+s)} r^{(\n-1)}_{r,s}(q^{2\m_1};q^2)
f_{\m+(r-s)\e_1}
\tag\eqname{\vgleltpeenr}
$$
and
$$
\p\bigl(R^{(\n+r+s)}_{l-r,m-s}(X_2,X_2^\ast,1;q^2)\bigr) f_{\m+(r-s)\e_1} =
 r^{(\n+r+s)}_{l-r,m-s}(q^{2\m_2};q^2) f_{\m+(l-m)\e_2+(s-r)\e}.
\tag\eqname{\vgleltptweer}
$$
Combination of \thetag{\vgleltpeenr} and \thetag{\vgleltptweer} gives
the explicit action of the first leg of the tensor product on the
right-hand side of the addition formula under the representation $\p$
on a basis vector.
\bigskip
We have now all ingredients to turn the  addition formula
\thetag{\vglabstraddform} into an addition formula in commuting variables
for $r_{l,m}^{(\n)}(\cdot;q)$.
Apply $\p\otimes\s$ to
\thetag{\vglabstraddform} and let the resulting
identity act on $f_\l\otimes e_{z+p-p\wedge\m_2}$ $(z\in \Zp)$ and take 
inner products
with $g_x(\m,p)$. Then use \thetag{\vglactionsltp}, \thetag{\vgleltpeenr},
\thetag{\vgleltptweer} and $\p\otimes\s$ being a $\ast$-representation
in combination with Corollary \thmref{\coractiolhs} to find
$$
\eqalign{
&(-q^{y+1})^{m-l} (q^{2N+2};q^2)_{m-l}^\hf (q^{2x-2l+2m+2};q^2)_{l-m}^\hf
\,  r_{l,m}^{(\n)}(q^{2x-2l+2m};q^2)\cr
&\qquad\qquad \times\ \langle f_\l\otimes e_{z+p-p\wedge\m_2},
g_{x+m-l}(\m+(m-l)\e_1,p)\rangle \cr
&=\sum_{r=0}^l\sum_{s=0}^m c_{l,m;r,s}^{(\n)}(q^2)\, (-q)^{r-s}
q^{(z+p-p\wedge\m_2+\l_2)(r+s)}\cr
&\qquad\qquad\times r^{(\n+r+s)}_{l-r,m-s}(q^{2\l_2};q^2) \, r_{r,s}^{(\n-1)}(q^{2\l_1};q^2) 
\, r^{(\n+r+s)}_{l-r,m-s}(q^{2(z+p-p\wedge\m_2)};q^2)\cr
&\qquad\qquad\times\langle f_{\l+(r-s)\e_1+(l-m+s-r)\e_2}\otimes
e_{z+p-p\wedge\m_2+l-m+s-r}, g_x(\m,p)\rangle
\cr}
\tag\eqname{\vgladdformeen}
$$
with $y=|p-\m_2|$, $N=\m_1+p\wedge\m_2$.
The inner products on both sides of \thetag{\vgladdformeen} can
be evaluated by use of Proposition \thmref{\propeigvectReen}. The left-hand
side is non-zero
if and only if $\l_2=z+\m_2 -p\wedge\m_2$, $\m_1+\m_2+m-l=\l_1+\l_2$ 
for some $z\in\Zp$ with $0\leq z\leq N+m-l$ and if so we find
$$
\langle f_\l\otimes e_{z+p-p\wedge\m_2}, g_{x+m-l}(\m+(m-l)\e_1,p)\rangle=
\hat K_z(q^{2x-2N};q^{2y},N+m-l;q^2).
\tag\eqname{\vglinprodeen}
$$
The inner product on the right-hand side of \thetag{\vgladdformeen} is non-zero
if and only if $\l_2=z+\m_2-p\wedge\m_2$, $\m_1+\m_2+m-l=\l_1+\l_2$
for some $z\in\Zp$ with
$0\leq z+l-m+s-r\leq N$ and if so we obtain
$$
\eqalign{
\langle f_{\l+(r-s)\e_1+(l-m+s-r)\e_2}\otimes e_{z+p-p\wedge\m_2+l-m+s-r},\
&g_x(\m,p)\rangle= \cr 
\hat K_{z+l-m+s-r}&(q^{2x-2N};q^{2y},N;q^2).
\cr}
\tag\eqname{\vglinprodtwee}
$$
If we interpret $\hat K_z(x;t,N;q)=0$ for $z<0$ and $z>N$, the inner products
in \thetag{\vglinprodeen} and \thetag{\vglinprodtwee} are well-defined
for all choices of $z\in\Z$.
We now choose $z,N$ such that the conditions $0\leq z\leq N+m-l$ and
$0\leq z+l-m+s-r\leq N$ are met for all possible choices of $r$ and $s$,
which is the case if we take $z+l\leq N$.

Plug \thetag{\vglinprodeen} and \thetag{\vglinprodtwee}, with the choices
$\l_2=z+\m_2 -p\wedge\m_2$ and $\m_1+\m_2+m-l=\l_1+\l_2$, in
\thetag{\vgladdformeen} to obtain an addition formula in commuting variables.
The dependence of the result on the variables $p, \m_1, \m_2$ is only on
$N= \m_1 + p\wedge\m_2$ and $y=|p-\m_2|$.
Observe that both sides are polynomial in $q^{2x}$ of degree 
$z+l$ which is at most equal to $N$, and 
that it holds for every $x\in\{0,\ldots,N\}$, so that
it still holds if we replace $q^{2x}$ by $x\in\C$. Next observe that if
we divide both sides by $q^{y(m-l-z)} \sqrt{(q^{2y+2};q^2)}$
we get an identity which is rational in $q^{2y}$ and holds for
$y\in\Zp$. Hence we can replace
$q^{2y}$ on both sides by $t$ with $0<t<q^{-2}$.
Finally, replacing the base $q^2$ by $q$ we obtain the addition formula in
commuting variables.

\proclaim{Theorem \theoremname{\thmaddqdisk}}
The functions $r^{(\n)}_{l,m}(x;q)$, $l,m\in\Zp$, $\n>0$,
defined by \thetag{\vglqdiskpolwithoutphase},
satisfy the addition formula for $x\in\C$, $t\in(0,q^{-1})$, $z,N\in\Zp$
with $z+l\leq N$,
$$
\multline
(-1)^{l-m}
\sqrt{(tq)^{m-l} {{(xq^{m-l+1};q)_{l-m}}\over{(q^{N+m-l+1};q)_{l-m}}} }
\  r_{l,m}^{(\n)}(xq^{m-l};q)\, \hat K_z(xq^{-N};t,N+m-l;q)\\
= \sum_{r=0}^l\sum_{s=0}^m c_{l,m;r,s}^{(\n)}(q)\, (-1)^{r-s}
q^{\hf (r-s)} t^{\hf(r+s)}q^{z(r+s)} \, r^{(\n+r+s)}_{l-r,m-s}(q^z;q)\,
r^{(\n+r+s)}_{l-r,m-s}(tq^z;q)\\
\times\, r_{r,s}^{(\n-1)}(q^{N+m-l-z};q)\,
\hat K_{z+l-m+s-r}(xq^{-N};t,N;q)
\endmultline
$$
with the notation as in \thetag{\vgldefortonormalaffqK} and
\thetag{\vgldefaffqK} and the constants $c_{l,m;r,s}^{(\n)}(q)$
given by \thetag{\cees}.
\endproclaim

\demo{Remark \theoremname{\remconstraintonz}}
If we drop the constraint $z+l\leq N$ the addition theorem remains valid
for $x$ replaced by $q^x \, (x=0,1,\hdots,N)$.
\enddemo

\demo{Remark \theoremname{\remthmaddformcommvar}} (i) The same result is obtained
if we work with the representation $\p\otimes\s_\phi$ of $\xt\otimes\yt$ for
arbitrary $\phi\in[0,2\p)$. Since $\p$ is a faithful $\ast$-representation
of $\xt$ and $\cap_\phi {\text{ker}}\,(\s_\phi$) is trivial, we see that
the addition formula in commuting variables is equivalent to
\thetag{\vglabstraddform}.

(ii) If we work with the representation $\p_1^\phi\otimes\s$
we obtain the limit case
$N\to\infty$ of Theorem \thmref{\thmaddqdisk} (cf. Remark
\thmref{\remthmpropeigvectReen}). Explicitly, we get
$$
\eqalign{
(-1)^{l-m} \sqrt{(tq)^{m-l} (xq^{m-l+1};q)_{l-m} }
\  &r_{l,m}^{(\n)}(xq^{m-l};q)\, w_z(xq^{m-l};t;q)\cr
= \sum_{r=0}^l\sum_{s=0}^m c_{l,m;r,s}^{(\n)}(q)\, (-1)^{r-s}
q^{\hf (r-s)} &t^{\hf(r+s)}q^{z(r+s)} \, r^{(\n+r+s)}_{l-r,m-s}(q^z;q)\cr
\times\, r^{(\n+r+s)}_{l-r,m-s}&(tq^z;q)\, w_{z+l-m+s-r}(x;t;q)
\cr}
\tag\eqname{\vgldegaddformWall}
$$
with the Wall polynomials defined by \thetag{\vgldefWallorthon}.
\enddemo

\demo{Remark \theoremname{\remaddformlittleqJ}} Taking $l=m$
in Theorem \thmref{\thmaddqdisk} gives an addition formula for
little $q$-Jacobi polynomials $p_l(\cdot;q^\n,1;q)$. In order to give 
it a similar form as the addition formula for the Jacobi
polynomial $P_l^{(\n,0)}$, cf. e.g. \cite{\VilKlim, Vol.~2, \S 11.4.2}, 
we split the double sum into two double sums;
$$
\aligned
&p_{l}^{(\n,0)}(x)\, \hat K_z(xq^{-N};t,N;q)=\\
&\sum_{r=0}^l\sum_{s=0}^r (1-\d_{rs}) c_{l,l;r,s}^{(\n)}(q)\, (-1)^{r-s}
q^{\hf (r-s)} t^{\hf(r+s)}q^{z(r+s)} 
(q^z,tq^z;q^{-1})_{r-s}^\hf (q^{N-z+1};q)_{r-s}^\hf \\
&\times p_{l-r}^{(\n+r+s,r-s)}(q^{z+s-r})
\,  p_{l-r}^{(\n+r+s,r-s)}(tq^{z+s-r})\,
p_s^{(\n-1,r-s)}(q^{N-z})\, 
\hat K_{z+s-r}(xq^{-N};t,N;q)\\
+&\sum_{r=0}^l\sum_{s=0}^r c_{l,l;r,s}^{(\n)}(q)\, (-1)^{r-s}
q^{\hf (r-s)} q^{r^2-s^2} t^{\hf(r+s)}q^{z(r+s)} 
(q^{z+1},tq^{z+1};q)_{r-s}^\hf (q^{N-z};q^{-1})_{r-s}^\hf \\
&\times p_{l-r}^{(\n+r+s,r-s)}(q^z)
\,  p_{l-r}^{(\n+r+s,r-s)}(tq^z)\, p_s^{(\n-1,r-s)}(q^{N-z+s-r})\, 
\hat K_{z+r-s}(xq^{-N};t,N;q),
\endaligned
%\tag\eqname{\vgladdformlittleqJ}
$$
where we use the notation $p^{(\a,\b)}_n(x)=p_n(x;q^\a,q^\b;q)$.
Koornwinder's \cite{\Koor} addition formula for the little $q$-Legendre
polynomials can be obtained from this identity
%\thetag{\vgladdformlittleqJ} 
by first
letting $N\to\infty$ and using the limit transition \thetag{\vgllimaffqKtoWall}
of the affine $q$-Krawtchouk polynomials to the Wall polynomials, or by
taking $l=m$ in the degenerate addition formula \thetag{\vgldegaddformWall}.
Next let $\n\to 0$. The double sums reduce to single sums, namely the sums
with $s=0$. This is because
$\lim_{\n\to 0} c_{l,m;r,s}^{(\n)}(q)$ is only non-zero if $r=0$ or $s=0$. With 
this observation Koornwinder's addition formula readily follows.
\enddemo

%%%%%%%%%%%%%%%%%%%%%%%%%%%%%%%%%%%%%%%%%%%%%%%%%%%
%%NEW SECTION %%%%%%%%%%%%%%%%%%%%%%%%%%%%%%%%%%%%%
%%%%%%%%%%%%%%%%%%%%%%%%%%%%%%%%%%%%%%%%%%%%%%%%%%%
\head\newsection . Product formula\endhead

In this section we derive a product formula from the addition formula
in Theorem \thmref{\thmaddqdisk}. For this we introduce a set of
orthogonal functions, which allow us to single out the $r=s=0$ term in the
addition formula.

For $a\in \Zp$ fixed and $s,r\in\Zp$  such that $0\leq a\leq N$
we define
$$
P_{r,s}^{(\nu)}(x,N;a,t;q) =
r_{r,s}^{(\nu)}(q^{N-a};q) \, \hat K_{a+s-r}(xq^{-N};t,N;q),
\tag\eqname{\vgldeforthofuncttwovar}
$$
where we use the notation as in \thetag{\vgldefortonormalaffqK} and
\thetag{\vglqdiskpolwithoutphase} and the convention that $\hat K_l (x;t,N;q)=0$
whenever $l<0$ or $l>N$.
Note that $P_{r,s}^{(\nu)}(x,N;a,t;q)$ is a polynomial in $x$.

\proclaim{Proposition \theoremname{\orthogonal}}
For $\n>-1$ the functions $P_{r,s}^{(\nu)}(x,N) = P_{r,s}^{(\nu)}(x,N;a,t;q)$
satisfy the orthogonality relations
$$
\sum_{N=a}^{\infty} \sum_{x=0}^N P_{r,s}^{(\nu)}(q^x,N)
P_{r\pr,s\pr}^{(\nu)}(q^x,N)
{{(q^N;q^{-1})_{N-x} (tq;q)_{N-x}}\over{(q;q)_{N-x} }} (tq)^x q^{(\n+1)(N-a)}
= \d_{rr\pr}\d_{ss\pr} h_{rs},
$$
where
$$
h_{rs} = {{(q;q)_r (q;q)_s}\over {(q^{\nu+1};q)_r (q^{\nu+1};q)_s}}
{{q^{(\nu+1)s}} \over {1-q^{\nu +r+s+1}}}.
$$
\endproclaim

\demo{Proof} This follows from the orthogonality relations
\thetag{\Kraworth} for the affine $q$-Krawtchouk polynomials,
the orthogonality relations for the little $q$-Jacobi polynomials,
\cite{\AndrA}, \cite{\GaspR};
$$
\sum_{x=0}^\infty {{(bq;q)_x}\over{(q;q)_x}}(aq)^x \big( p_lp_m\bigr)
(q^x;a,b;q) =\d_{lm}
{{(q,bq;q)_l (aq)^l (1-abq) (abq^2;q)_\infty}\over
{(aq,abq;q)_l (1-abq^{2l+1}) (aq;q)_\infty}},
$$
with $0<aq<1$, $b<q^{-1}$ and the definition
\thetag{\vglqdiskpolwithoutphase}. \qed\enddemo

The addition formula of Theorem \thmref{\thmaddqdisk} can be viewed as
the development of the left-hand side in terms of the orthogonal functions
$P_{r,s}^{(\nu-1)}(x,N;z+l-m,t;q)$. Using the orthogonality relations
of Proposition \thmref{\orthogonal} we can pick out the term $r=s=0$ to find
a product formula.

\proclaim{Theorem \theoremname{\thmprodformlqdiskpol}}
The functions $r^{(\n)}_{l,m}(x;q)$, $l,m\in\Zp$, $\n>0$,
defined by \thetag{\vglqdiskpolwithoutphase},
satisfy the following product formula for
$t\in(0,q^{-1})$, $z\in\Zp$;
$$
\eqalign{
r^{(\n)}_{l,m}(q^z;q)\, r^{(\n)}_{l,m}(tq^z;q) = (&1-q^{\n})\ 
\sum_{N=z+l-m}^\infty \sum_{x=0}^N (-1)^{l-m}
\sqrt{(tq)^{m-l} {{(q^{x-l+m+1};q)_{l-m}}\over{(q^{N+m-l+1};q)_{l-m}}} }
\cr
&\times {{(q^N;q^{-1})_{N-x} (tq;q)_{N-x}}\over{(q;q)_{N-x} }}
(tq)^x q^{\n(N-z-l+m)}\,  r_{l,m}^{(\n)}(q^{x+m-l};q) \cr 
&\times
\hat K_z(q^{x-N};t,N+m-l;q)\, \hat K_{z+l-m}(q^{x-N};t,N;q)
\cr}\tag\eqname{\prodformlidentity}
$$
with the notation as in \thetag{\vgldefortonormalaffqK}.
\endproclaim

Notice that the inner summation in fact runs from the maximum of $0$ and $l-m$
to $N$.

\demo{Remark \theoremname{\remonsymmetryprodforml}}
The symmetry relation $r^{(\n)}_{l,m}(xq^{m-l};q) = r^{(\n)}_{m,l}(x;q)$
gives rise to a symmetry of the product formula \thetag{\prodformlidentity}.
Namely, changing $x,z,N$ to $x+l-m, z+m-l, N+l-m$ on both sides of 
\thetag{\prodformlidentity} has the same effect as interchanging the
parameters $l$ and $m$. In particular the full identity 
\thetag{\prodformlidentity} is known when it is known either for
$l\leq m$ or $l\geq m$.
\enddemo

By specialising $l=m$ we obtain the following corollary for the little
$q$-Jacobi polynomials.

\proclaim{Corollary \theoremname{\thmprodformllittleqJacobi}}
The little $q$-Jacobi polynomials $p_l(x;q^\n,1;q)$, $l\in\Zp$,
$\n>0$, defined by \thetag{\vgldeflittleqJ},
satisfy the following product formula for
$t\in(0,q^{-1})$, $z\in\Zp$;
$$
\multline
p_l(q^z;q^\n,1;q)\, p_l(tq^z;q^\n,1;q) = (1-q^{\n})\ 
\sum_{N=z}^\infty \sum_{x=0}^N
{{(q^N;q^{-1})_{N-x} (tq;q)_{N-x}}\over{(q;q)_{N-x} }} (tq)^x q^{\n(N-z)}
\\ \times p_l(q^x;q^\n,1;q)\, \Bigl(\hat K_z(q^{x-N};t,N;q)\Bigr)^2.
\endmultline
$$
\endproclaim

\demo{Remark \theoremname{\remprodformlittleqJ}} 
One can show, by use of Abel's partial summation formula, 
that in the limit $\nu\to 0$ the right-hand side in 
Corollary \thmref{\thmprodformllittleqJacobi} tends to
$$
\sum_{x=0}^{\infty} \frac{(q^{x+1};q)_{\infty} (tq;q)_{\infty}}
{(q;q)_{\infty}} (tq)^x  p_l(q^x;1,1;q) \Bigl(w_z(q^x;t;q\Bigr)^2
$$
with $w_z(\cdot;t;q)$ as in \thetag{\vgldefWallorthon}. After some
rewriting it is easily seen that for $t=q^y$ this is the same 
expression as the one derived by Koornwinder \cite{\Koor, Thm.~5.1}.
\enddemo

%%%%%%%%%%%%%%%%%%%%%%%%%%%%%%%%%%%%%%%%%%%%%%%%%%%
%%NEW SECTION %%%%%%%%%%%%%%%%%%%%%%%%%%%%%%%%%%%%%
%%%%%%%%%%%%%%%%%%%%%%%%%%%%%%%%%%%%%%%%%%%%%%%%%%%
\head\newsection . The limit case $q\uparrow 1$\endhead

In this section we consider the limit case $q\uparrow 1$ of the addition
and product formula for the little $q$-disk polynomials to the
addition and product formula \thetag{\classadd} and \thetag{\classprod}
for the disk polynomials.
For this limit transition we use the methods developed by Van Assche and
Koornwinder \cite{\AsKoo}. Since some of the calculations are rather
tedious we restrict ourselves in this section to giving the key steps; 
details can be found in the Appendix.
Calculations similar to the ones appearing here were made 
in \cite{\KoelSwar}.

First we apply \cite{\AsKoo, Thm.~1} to the polynomial
$\hat K_l(xq^{-N};t,N+a;q)$. It gives the limit of the quotient of two
such polynomials as $q$ tends to $1$. Although \cite{\AsKoo, Thm.~1} is not 
formulated for discrete orthogonal polynomials, 
the proof of \cite{\AsKoo, Thm.~1} is also applicable to discrete orthogonal 
polynomials provided that $N$, and hence the number of orthogonal polynomials, 
tends to infinity.
It is straightforward to verify
that the conditions of \cite{\AsKoo, Thm.~1} are satisfied.
Explicitly, we have the following result.

\proclaim{Proposition \theoremname{\proplimaffqKrawt}} For $0<t<c^{-1}$,
$c\in (0,1)$, $a,m\in\Z$ and $N\in \N$, we have
$$
\lim_{n\to\infty} {{\hat K_{n+m}(x(c^{1/n})^{-nN};t,nN+a;c^{1/n})}\over
{\hat K_n(x(c^{1/n})^{-nN};t,nN+a;c^{1/n})}} =
\rho^m \bigl( (x-B)/(2A)\bigr),
$$
uniformly for $x$ on compact subsets of $\C \backslash[0,c^{-a}]$. Here
$$
A = c\sqrt{t(1-tc)(1-c)(1-c^{N-1})}>0, \qquad
B = c + tc -2tc^2 + tc^{1+N}\in\R,
$$
and $\rho(x) = x+\sqrt{x^2-1}$ with the square root taken such that
$|\rho(x)|>1$ for $x\in\C\backslash[-1,1]$.
\endproclaim

In order to apply Proposition \thmref{\proplimaffqKrawt} to the addition
formula we need the connection coefficients between affine $q$-Krawtchouk
polynomials for different $N$.

\proclaim {Proposition \theoremname{\connection}}
For $a\in\Z$, $a\leq N$, $x\in\C$ and $0\leq l\leq N\wedge (N-a)$ we have the
connection formula
$$
K_l (x;t,N;q) = \sum_{k=0}^l {{(q^l;q^{-1})_k}\over{(q;q)_k}}
q^{(-N+a)k} \frac{(q^{-a};q)_k}{(q^{-N};q)_l} (q^{-N+a};q)_{l-k}
K_{l-k} (x;t,N-a;q).
$$
\endproclaim

\demo{Proof} There are several ways of proving the proposition. We use the
generating function for the affine $q$-Krawtchouk polynomials; for
$x\in\{ 0,\ldots,N\}$,
$$
(zq^{-N};q)_{N-x} {}_{1}\vp_1(q^{-x};tq;q,tqz)=
\sum_{l=0}^N \frac{(q^{-N};q)_l}{(q;q)_l} K_l(q^{-x};t,N;q) z^l.
\tag\eqname{\Krawgenfunction}
$$
This identity is proved as follows: insert \thetag{\vgldefaffqK} in the 
right-hand side of \thetag{\Krawgenfunction}, change summations and apply
the $q$-binomial theorem
to find that this right-hand side equals
$$
(zq^{-N};q)_N\, {}_3\phi_2 \left [ {q^{-x}, 0, 0 \atop
tq , q/z}; q,q \right ].
$$
Then use the limiting case $b\to\infty$ of the identity which one gets by equating
\cite{\GaspR, (III.5)} and the right-hand side of \cite{\GaspR, (III.4)}, 
with $a=q^{-x}, c=tq$ and $z$ replaced by $tqz/b$, to obtain the left-hand 
side of \thetag{\Krawgenfunction}. Now, for $|z|<q^{N-a}$
$$
\eqalign{
&\sum_{l=0}^{N} \frac{(q^{-N};q)_l}{(q;q)_l} K_l(q^{-x};t,N;q) z^l =\cr
&(zq^{-N};q)_{a} (zq^{-N+a};q)_{N-a-x}
{}_{1}\vp_1(q^{-x};tq;q,tqz) =\cr
&\sum_{k=0}^{\infty} \frac{(q^{-a};q)_k}{(q;q)_k} (zq^{-N+a})^k\
\sum_{j=0}^{N-a} \frac{(q^{-N+a};q)_j}{(q;q)_j}  K_j(q^{-x};t,N-a;q) z^j= \cr
&\sum_{l=0}^{N} \bigl( \sum_{k=0}^l \frac{(q^{-a};q)_k}{(q;q)_k}
q^{(-N+a)k} \frac{(q^{-N+a};q)_{l-k}}{(q;q)_{l-k}}  K_{l-k}(q^{-x};t,N-a;q)
\bigr) z^l
\cr}
$$
by the $q$-binomial theorem and \thetag{\Krawgenfunction}. Observe that
the right-hand side is in fact a polynomial because of the convention that 
$K_l (x;t,N;q) =0$ whenever $l>N$. This proves the
proposition for $x$ replaced by $q^{-x}$, 
$x\in\{0,\ldots,N\wedge (N-a)\}$. Since both
sides are polynomials of degree not exceeding $N\wedge (N-a)$ the result
follows for arbitrary $x\in\C$. \qed\enddemo

Next use Propositions \thmref{\proplimaffqKrawt} and \thmref{\connection}
and the binomial theorem
to obtain the limit of the quotient of the affine $q$-Krawtchouk polynomials
appearing in the addition formula of Theorem \thmref{\thmaddqdisk}.
We replace $q$, $z$, $N$ by $c^{1/nz}$, $nz$, $nzN$ for $c\in(0,1)$ and
$z,N,n \in\N$ and apply the binomial theorem to find
$$
\multline
\lim_{n\to\infty} {{\hat
K_{nz+l-m+s-r}(x(c^{1/nz})^{-nzN};t,nzN;c^{1/nz})}\over
{\hat K_{nz}(x(c^{1/nz})^{-nzN};t,nzN+m-l;c^{1/nz})}} = \\
\bigl( c \sqrt{{1-c^{N-1}}\over {1-c^N}}\bigr)^{l-m} 
\rho^{l-m+s-r}\bigl( {{x-B}\over {2A}}\bigr)
\biggl(  1+ {{\sqrt{(1-tc)(1-c)(1-c^{N-1})}\over{\sqrt{t} c^N
(1-c^{-N+1}) \rho({{x-B}\over {2A}})}}}\biggr) ^{l-m}
\endmultline
$$
uniformly in $x$ on compact sets of $\C\backslash [0, c^{m-l}]$.
Here $A$ and $B$ are as in 
Proposition \thmref\proplimaffqKrawt. We remark that in case $m-l\geq 0$ 
the sum obtained on the left from applying
the connection formula of Proposition \thmref{\connection} to the numerator 
polynomial is finite, even when we let
$n$ tend to infinity. If $m-l\leq 0$, one develops the denominator
polynomial as a (finite) sum by means of the connection formula. 

Introduce the parameter $u>1$ by putting
$(x-B)/ (2A) = {\hf} (u+ u^{-1})$. Then
$\rho((x-B)/(2A)) =u$. If we now divide both sides of
the addition formula in Theorem \thmref{\thmaddqdisk} by the polynomial 
$\hat K_z(xq^{-N};t, N+m-l;q))$
and, as before, replace
$q$, $z$, $N$ by $c^{1/nz}$, $nz$, $nzN$
and let $n$ tend to infinity, then, after some careful computations
using the above result and an analytic continuation argument replacing 
$- u^{-1}$ by 
$e^{i\psi}$, we end up with the identity  
$$
\eqalign{
R_{l,m}^{(\nu)}(z) = 
\sum_{r=0}^l\sum_{s=0}^m C_{l,m;r,s}^{(\nu)} &(\sin\theta_1)^{r+s}
 R_{l-r,m-s}^{(\nu +r+s)}(\cos\theta_1)\cr
\times &(\sin\theta_2)^{r+s}  R_{l-r,m-s}^{(\nu +r+s)}(\cos\theta_2)
 R_{r,s}^{(\nu -1)}(\rho e^{i\psi}).
\cr}
\tag\eqname{\classaddpartial}
$$
Here 
$\cos\theta_1=\sqrt{1-c}$, $\cos\theta_2=\sqrt{1-tc}$,
$\rho=\sqrt{1-c^{N-1}}$ and 
$ z= \sqrt{1-c} \sqrt{1-tc} + c\sqrt{t}\sqrt{1-c^{N-1}} e^{i\psi}=
\cos\theta_1\cos\theta_2 + \sin\theta_1\sin\theta_2\rho e^{i\psi}$.
In order to obtain the general classical addition formula for disk polynomials
\thetag{\classadd} from this, one simply observes that multiplication 
of both sides
of \thetag{\classaddpartial} with $e^{i(l-m)\phi_1} e^{i(l-m)\phi_2}$ by
homogeneity of the polynomials has the effect of replacing $\cos\theta_1$ by 
$e^{i\phi_1}\cos\theta_1$ and
$\cos\theta_2$ by $e^{i\phi_2}\cos\theta_2$, and hence produces the full
classical identity.

\medskip

Let us turn now to the product formula. In view of Remark 
\thmref{\remonsymmetryprodforml} we can restrict ourselves to the case $l\leq m$.
Observe that the right-hand side of \thetag{\prodformlidentity}
is of the form
$$
(1-q^{\n}) \sum_{N=z+l-m}^{\infty}\, q^{\n (N-z)}\, g(N,z,q)
\tag\eqname{\shorthand}
$$
with
$$
\eqalign{
g(N,z,q) = \sum_{x=0}^N &(-1)^{l-m}
\sqrt{(tq)^{m-l} {{(q^{x-l+m+1};q)_{l-m}}\over{(q^{N+m-l+1};q)_{l-m}}} }
\cr
&\times {{(q^N;q^{-1})_{N-x} (tq;q)_{N-x}}\over{(q;q)_{N-x} }}
(tq)^x \,  r_{l,m}^{(\n)}(q^{x+m-l};q) \cr 
&\times
\hat K_z(q^{x-N};t,N+m-l;q)\, \hat K_{z+l-m}(q^{x-N};t,N;q).\cr}
\tag\eqname{\functionsg}
$$
If we replace $q,z,N$ in \thetag{\shorthand} by $c^{1/nz}, nz, nzN$, 
with $c\in (0,1)$, and let
$n$ tend to infinity, we formally get 
$$
\eqalign{
&\lim_{n\to\infty} nz\, (1-c^{\nu /nz}) \sum_{{N=1}\atop 
{\text{step-size } 1/nz}}^{\infty}
c^{\n (N-1)}\, \frac{g(nzN, nz, c^{1/ nz})}{nz} =\cr
&-\n \log c \int_{N=1}^{\infty} c^{\n (N-1)}\, 
\lim_{n\to\infty}\bigl( \frac{g(nzN, nz, c^{1/ nz})}{nz}\bigr)\, dN =\cr}
$$
$$
\eqalign{
&-\n \log c \int_{N=1}^{\infty} c^{\n (N-1)}\, G(\sqrt{1-c^{N-1}})\, dN=\cr
&\nu\, \int_0^1 G(\sqrt{1-\r})\, \r^{\n -1}\, d\r\cr}
\tag\eqname{\formalresult}
$$
where we made the substitution $\r = c^{N-1}$ in the last step. 
To calculate the limit function $G$ and to show that it, as a function
of $N$, in fact only depends on $\sqrt{1-c^{N-1}}$, we proceed 
as follows.
Use the connection formula of Proposition \thmref{\connection} to replace 
the Krawtchouk
polynomial $\hat K_z (q^{x-N};t,N+m-l;q)$ in \thetag{\functionsg}
by a sum of polynomials 
$\hat K_{z-k} (q^{x-N};t,N;q)$ $(k=0,\hdots ,z)$.
Then replace $q,z,N$ by $c^{1/nz}$, $nz$, 
$nzN$ for $c\in(0,1)$ and send $n$ to infinity. According to 
\cite{\AsKoo, Thm.~2}, which again is valid also for discrete orthogonal
polynomials, and the  binomial theorem we get 
$$
\eqalign{
G(N) = (-&1)^{l-m}\, [c\sqrt{t (1-c^{N-1})}]^{m-l} \sum_{k=0}^{m-l}
(-1)^k t^{-\frac{k}{2}} c^{-Nk} \frac{(l-m)_k}{k!}\cr
&\times [(1-c)(1-tc)(1-c^{N-1})]^{\frac{k}{2}} (1-c^{-N+1})^{-k}\cr
&\times \frac{1}{\pi} \int_{B-2A}^{B+2A} (1-x)^{\frac{l-m}{2}}
R_{l,m}^{(\nu)} (\sqrt{1-x}) \frac{T_{m-l-k}\bigl( \frac{x-B}{2A}\bigr)}
{\sqrt{4A^2 - (x-B)^2}} dx\cr}
$$
where $T_n$ are Chebyshev polynomials of the first kind and $A,B$ are as 
in Proposition 
\thmref{\proplimaffqKrawt}.
Next perform the substitution $(x-B)/(2A) = \hf (u+u^{-1}) =
\cos\psi$ with $u= e^{i\psi}$. Then interchange the order of summation and 
integration and use that $T_n(\cos\psi) = \hf (e^{in\psi} + e^{-in\psi})$
to get two finite sums. 
Apply the binomial theorem to both of these sums to obtain
$$
\eqalign{
G(N) =\, &(-1)^{l-m}\, [c\sqrt{t (1-c^{N-1})}]^{m-l} \cr
&\times \frac{1}{2\pi} \int_{0}^{\pi} 
(1-(2A\cos\psi +B))^{\frac{l-m}{2}}
R_{l,m}^{(\nu)} (\sqrt{1-(2A\cos\psi +B)}) \cr
&\times \biggl\{ \biggl(e^{i\psi} + \sqrt{t^{-1}(1-c)(1-tc)(1-c^{N-1})}
\frac{c^{-N}}{1-c^{-N+1}}\biggr)^{m-l} \cr
&\qquad + \biggl(e^{-i\psi} + \sqrt{t^{-1}(1-c)(1-tc)(1-c^{N-1})}
\frac{c^{-N}}{1-c^{-N+1}}\biggr)^{m-l}\biggr\} d\psi.\cr}
\tag\eqname{\firststep}
$$
As we assumed that $l\leq m$, we have
$$
(1-(2A\cos\psi +B))^{\frac{l-m}{2}}
R_{l,m}^{(\nu)} (\sqrt{1-(2A\cos\psi +B)}) =
 P_l^{(\nu,l-m)} (1-2(2A\cos\psi +B))
$$
with $A,B$ as in Proposition
\thmref{\proplimaffqKrawt}. It is easily checked that the following equality
holds;
$$
\eqalign{
1-(2A\cos\psi +B) = &\bigl( \sqrt{(1-c)(1-tc)} - c\sqrt{t(1-c^{N-1})}
e^{i\psi}\bigr) \cr
&\times \bigl( \sqrt{(1-c)(1-tc)} - c\sqrt{t(1-c^{N-1})}
e^{-i\psi}\bigr). \cr}
$$
When we use this we find that the right-hand side of \thetag{\firststep} 
reduces to
$$
\eqalign{
\frac{1}{2\pi} \int_{0}^{\pi} P_{l\wedge m}^{(\nu, |l-m|)}
(1-2(2A\cos\psi &+B)) 
\biggl\{ \biggl( \sqrt{(1-c)(1-tc)} - c\sqrt{t(1-c^{N-1})}
e^{-i\psi}\biggr) ^{m-l}\cr
&+ \biggl( \sqrt{(1-c)(1-tc)} - c\sqrt{t(1-c^{N-1})}
e^{i\psi}\biggr) ^{m-l} \biggr\} d\psi = \cr 
\frac{1}{2\pi} \int_{0}^{2\pi} P_{l\wedge m}^{(\nu, |l-m|)}
(1-2(2A\cos\psi &+B))  \biggl( \sqrt{(1-c)(1-tc)} + c\sqrt{t(1-c^{N-1})}
e^{-i\psi}\biggr) ^{m-l} d\psi \cr}
$$
where in the second step we combined the two integrals and changed 
$\psi$ to $\psi -\pi$.\par
So indeed we see that $G$, as a function of $N$, only depends on
$\sqrt{1-c^{N-1}}$. Combining this result with \thetag{\formalresult} we find 
that we end up with
$$
\eqalign{
\frac{\nu}{2\pi} \int_0^1 \int_{0}^{2\pi}
&P_{l\wedge m}^{(\nu, |l-m|)}(1-2(2A\cos\psi +B))\cr
&\times \biggl( \sqrt{(1-c)(1-tc)} + c\sqrt{t}\sqrt{1-\rho}
e^{-i\psi}\biggr) ^{m-l} \rho^{\nu-1} d\psi d\rho =\cr
\frac{\nu}{\pi} \int_0^1 \int_{0}^{2\pi}
&P_{l\wedge m}^{(\nu, |l-m|)}(2 |z|^2 -1) (\bar z)^{m-l} r (1-r^2)^{\nu -1}
d\psi dr=\cr
\frac{\nu}{\pi} \int_0^1 \int_{0}^{2\pi}
&R_{l,m}^{(\nu)} (z) r  (1-r^2)^{\nu -1} d\psi dr\cr}\tag\eqname{\righthand}
$$
where we substituted $r=\sqrt{1-\r}$ and where $z= \sqrt{(1-c)(1-tc)} 
+ c\,\sqrt{t}\, r\, e^{i\psi}$. 

It is easily seen from the limit transition of the little 
$q$-Jacobi 
polynomials to the Jacobi polynomials that, under the above substitutions, 
the left-hand side of the product formula \thetag{\prodformlidentity}
tends to the product
$$
R_{l,m}^{(\nu)}(\sqrt{1-c})\, R_{l,m}^{(\nu)}(\sqrt{1-tc})\tag\eqname{\lefthand}
$$ 
as $n$ tends to infinity. From the fact that \thetag{\righthand} equals
\thetag{\lefthand} we obtain the full classical product
formula \thetag{\classprod} when we substitute
$\cos\theta_1$ and $\cos\theta_2$
for $\sqrt{1-c}$ and $\sqrt{1-tc}$ respectively, and invoke the same homogeneity
argument as in the case of the limit transition of the addition formula.
Since the case $l\geq m$ followed from a symmetry argument, we are done.

%\vfill\eject
\bigskip

%%%%%%%%%%%%%%%%%%%%%%%%%%%%%%%%%%%%%%%%%%%%%%%%%
%%%%%%%%%%%%%%%%%%%%%%%%%%%%%%%%%%%%%%%%%%%%%%%%%
%%%%%% APPENDIX
%%%%%%%%%%%%%%%%%%%%%%%%%%%%%%%%%%%%%%%%%%%%%%%%%
%%%%%%%%%%%%%%%%%%%%%%%%%%%%%%%%%%%%%%%%%%%%%%%%%

\def\newsection{\global\advance\sectionno by 1
                \global\eqnumber=1
                \global\theoremno=1
                \global\cntsubsecno=0}

\def\eqname#1{A.\the\eqnumber
              \xdef#1{{A.\the\eqnumber}}
              \global\advance\eqnumber by 1}

\head\newsection  Appendix\endhead
\medskip
In this appendix we give justifications for, and details of, the calculations
that were made in section 7.\par
To prove Proposition \thmref{\proplimaffqKrawt} we have to show that
all the conditions of \cite{\AsKoo, Thm.~1} are satisfied. From the 
recurrence relation \thetag{\recurrorthonormKraw} we obtain
$(a\in\Z,\, z,N,n\in\N,\, c\in (0,1))$
$$
\eqalign{
&x\hat K_{nz}(xc^{-N};t,nzN+a;c^{1/nz}) = \cr
&- c^N a_{nz}(t,nzN+a;c^{1/nz}) \hat K_{nz+1}(xc^{-N};t,nzN+a;c^{1/nz}) \cr
&+ c^N\bigl( 1- b_{nz}(t,nzN+a;c^{1/nz})\bigr) 
\hat K_{nz}(xc^{-N};t,nzN+a;c^{1/nz}) \cr
&- c^N a_{nz-1}(t,nzN+a;c^{1/nz}) \hat K_{nz-1}(xc^{-N};t,nzN+a;c^{1/nz})\cr}
$$
with $a_{l}(t,N;q)$ and $b_{l}(t,N;q)$ given by \thetag{\recurrcoeff}. Put
$$
\hat a_l(t,N;q) = -q^N  a_l(t,N;q)\, ,\qquad
\hat b_l(t,N;q) = q^N (1-b_l(t,N;q)).
$$
Then it is immediate that
$$
\eqalign{
&\lim_{n\to\infty} \hat a_{nz}(t,nzN+a;c^{1/nz}) = A\cr
&\lim_{n\to\infty} \hat b_{nz}(t,nzN+a;c^{1/nz}) = B\cr}
$$
with
$$
A = c\sqrt{t(1-tc)(1-c)(1-c^{N-1})}>0, \qquad
B = c + tc -2tc^2 + tc^{1+N}\in\R\, ;
$$
see Proposition \thmref{\proplimaffqKrawt}. Furthermore, the expressions
$ \hat a_{l} (t,nzN+a;c^{1/nz})^2$ and $ \hat b_{l} (t,nzN+a;c^{1/nz})$
are both of the form $f(c^{1/nz}, c^{l/nz})$ with $f$ a polynomial.
Hence, by Lipschitz continuity, we find
that
$$
\eqalign{
&\lim_{n\to\infty} \bigl[ \hat a_{l}(t,nzN+a;c^{1/nz})^2 -
\hat a_{l-1}(t,nzN+a;c^{1/nz})^2\bigr] = 0\cr
&\lim_{n\to\infty} \bigl[ \hat b_{l}(t,nzN+a;c^{1/nz}) -
\hat b_{l-1}(t,nzN+a;c^{1/nz})\bigr] =0\cr}
$$
uniformly in $l$. The support of the orthogonality measure of the polynomials
$\hat K_{l}(x;t,N;q)$ is contained in the interval $[1,q^{-N}]$, hence the
polynomials $x\mapsto \hat K_{l}(xc^{-N};t,nzN+a;c^{1/nz})$ have a measure
of orthogonality
with support in $[c^N, c^{-a/nz}]$. A given compact subset of
$\C\bs [0,1]$ will have no intersection with $[c^N, c^{-a/nz}]$ if $n$ is big
enough (this is trivial when $a\leq 0$). From the proof of \cite{\AsKoo, Thm.~1} 
it follows that the statement of that theorem is still valid for orthogonal
polynomials with finite support, provided that the support of the
orthogonality measure for the polynomials $p_l(\cdot;n)$ contains more 
than $n+2$ points, which is true in our case. In this way we have verified 
that all the
conditions of \cite{\AsKoo, Thm.~1} are fulfilled, which proves Proposition
\thmref{\proplimaffqKrawt}.\par
As to the explicit calculations, observe that from Proposition
\thmref{\connection} we obtain the following connection formula for the
orthonormal affine $q$-Krawtchouk polynomials;
$$
\eqalign{
\hat K_{l}(x;t,N;q) = \sum_{k=0}^l &(-1)^k (tq)^{-\frac{k}{2}} q^{(-N+a)k}
\frac{(q^{-a};q)_k (q^{-N+a};q)_{l-k}}{(q;q)_k (q^{-N};q)_l}\cr
& \times \bigl[ (tq^l;q^{-1})_k (q^N;q^{-1})_a (q^{N-l+1};q)_{k-a} 
(q^l;q^{-1})_k\bigr]^{\hf}\cr
& \times \hat K_{l-k}(x;t,N-a;q)\cr
= \sum_{k=0}^l &(-1)^k (tq)^{-\frac{k}{2}} q^{(-N+a)k}
\frac{(q^{-a};q)_k}{(q;q)_k (q^{-N};q)_a (q^{-N+l-1};q^{-1})_{k-a}}\cr
&\times\bigl[ (tq^l;q^{-1})_k (q^N;q^{-1})_a (q^{N-l+1};q)_{k-a} 
(q^l;q^{-1})_k\bigr]^{\hf}\cr
& \times\hat K_{l-k}(x;t,N-a;q)\cr}\tag\eqname{\connformorthonormKraw}
$$
for $a\leq N,\, 0\leq l\leq N\wedge (N-a)$. Now consider the quotient
$$
\frac{\hat K_{z+l-m+s-r}(xq^{-N};t,N;q)}
{\hat K_{z}(xq^{-N};t,N+m-l;q)}.\tag\eqname{\quotientKraw}
$$
Suppose $a=l-m\geq 0$. Then the sum in \thetag{\connformorthonormKraw} will be
a finite sum $k=0,\hdots,a$ because of the factor $(q^{-a};q)_k$. So if we
develop the numerator of \thetag{\quotientKraw} in terms of polynomials
$\hat K_{z+l-m+s-r-k}(xq^{-N};t,N;q)$ $(k=0,\hdots, z+l-m+s-r)$ this sum
will always be finite in this case. Replacing $q,z,N$ by $c^{1/nz}, nz, nzN$
with $c\in (0,1)$ and $n,z,N\in\N$ and applying Proposition
\thmref{\proplimaffqKrawt} we then obtain;
$$
\eqalign{
\lim_{n\to\infty} &\frac{\hat
K_{nz+l-m+s-r}(x(c^{1/nz})^{-nzN};t,nzN;c^{1/nz})}
{\hat K_{nz}(x(c^{1/nz})^{-nzN};t,nzN+m-l;c^{1/nz})} =\cr
\sum_{k=0}^{l-m} &(-1)^k t^{-\hf k} c^{-Nk} \frac{(m-l)_k}{k!}
(1-c^{-N})^{m-l} (1-c^{-N+1})^{l-m-k}\cr
&\times\bigl[ (1-tc)^k (1-c^N)^{l-m} (1-c^{N-1})^{k-l+m} (1-c)^k\bigr]^{\hf}
\r^{l-m+s-r-k}\bigl( \frac{x-B}{2A}\bigr)=\cr
&\!\!\!\!\!\!\!\!\!\! \bigr[ \frac{1-c^{-N+1}}{1-c^{-N}}\bigr] ^{l-m} 
\bigr[ \frac{1-c^{N}}{1-c^{N-1}}\bigr] ^{\hf (l-m)} 
\r^{l-m+s-r}\bigl( \frac{x-B}{2A}\bigr)\cr
&\times\sum_{k=0}^{l-m}  \frac{(m-l)_k}{k!} \bigl[ - 
\frac{\sqrt{t^{-1}(1-c)(1-tc)
(1-c^{N-1})}}{c^N (1-c^{-N+1}) \r (\frac{x-B}{2A})}\bigr]^k =\cr
c^{l-m} &\bigr[ \frac{1-c^{N-1}}{1-c^{N}}\bigr] ^{\hf (l-m)} 
\r^{l-m+s-r}\bigl( \frac{x-B}{2A}\bigr)\cr
&\times\bigl[ 1+ \frac{\sqrt{t^{-1}(1-c)(1-tc)
(1-c^{N-1})}}{c^N (1-c^{-N+1}) \r (\frac{x-B}{2A})}\bigr]^{l-m} 
\cr}\tag\eqname{\limitquotients}
$$
where in the last step we applied the binomial theorem and where $A,B$ are as in
Proposition \thmref{\proplimaffqKrawt}. The case $a=l-m\leq 0$ is treated
similarly; in that case develop the numerator of \thetag{\quotientKraw} using the
connection formula \thetag{\connformorthonormKraw} (so in fact one does the
same procedure as above for the inverse of the quotient 
\thetag{\quotientKraw}).
Let us assume from now on that $a=l-m\geq 0$. An easy calculation shows that
$$
\lim_{q\uparrow 1} r_{l,m}^{(\nu)} (x;q) = R_{l,m}^{(\nu)}(\sqrt{1-x})
$$
with the disk polynomial $R_{l,m}^{(\nu)}(z)$ as in \thetag{\defdiskpol}.
If we now divide both sides of the addition theorem \thmref{\thmaddqdisk} by
the polynomial $\hat K_z(xq^{-N};t,N+m-l;q)$ and make the substitutions
for $q,z,N$ as before, the left-hand side will tend to
$$
(-1)^{l-m} \bigl[ \frac{1-x}{t(1-c^N)}\bigr] ^{\hf (l-m)} 
R_{l,m}^{(\nu)}(\sqrt{1-x})\tag\eqname{\LHside}
$$
as $n\to\infty$. In the same way, using \thetag{\limitquotients}, the right-hand
side tends to
$$
\eqalign{
\sum_{r=0}^l \sum_{s=0}^m &C_{l,m;r,s}^{(\n)}\, (-1)^{r-s} t^{\hf (r+s)} c^{r+s}\,
R_{l-r,m-s}^{(\nu +r+s)}(\sqrt{1-tc})\cr
&\times R_{l-r,m-s}^{(\nu +r+s)}(\sqrt{1-c})\,
R_{r,s}^{(\nu -1)}(\sqrt{1-c^{N-1}})\,
 c^{l-m}\, \bigl[ \frac{1-c^{N-1}}{1-c^N}\bigr] ^{\hf (l-m)} \cr
&\times\r^{l-m+s-r} (\frac{x-B}{2A})\, \bigl[ 1+ \frac{\sqrt{t^{-1}(1-c)(1-tc)
(1-c^{N-1})}}{c^N (1-c^{-N+1}) \r (\frac{x-B}{2A})}\bigr]^{l-m} \cr}
\tag\eqname{\RHside}
$$
with $C_{l,m;r,s}^{(\n)}$ as in \thetag{\classadd}. Now put $(x-B)/(2A) =
\hf (u+u^{-1})$, so that $\r ((x-B)/(2A)) = u$. Write
$$
F(\psi; t,N;c) = c^{l-m} \bigl[ \frac{1-c^{N-1}}{1-c^N}\bigr] ^{\hf (l-m)} 
\bigl[ 1 + \frac{\sqrt{t^{-1}(1-c)(1-tc)
(1-c^{N-1})}}{c^N (1-c^{-N+1}) u} \bigr]^{l-m}.
$$
It is a straightforward exercise to verify that
$$
\eqalign{
&(-1)^{l-m} \bigl[ \frac{1-x}{t(1-c^N)}\bigr] ^{\hf (l-m)} F(\psi; t,N;c)^{-1} =
\cr
&(1-x)^{\hf (m-l)} (\sqrt{(1-c)(1-tc)} - u^{-1} c \sqrt{t(1-c^{N-1})})^{l-m}.
\cr}\tag\eqname{\tussenstapje}
$$
Moreover, from $ (x-B)/(2A)=\hf (u+u^{-1})$ and from the values of $A$ and $B$ it
follows that, when $z= \sqrt{(1-c)(1-tc)} - u^{-1} c \sqrt{t(1-c^{N-1})}$
then $|z|=\sqrt{1-x}$. Now we can put all the pieces 
together: equate
\thetag{\LHside} and \thetag{\RHside} and divide both sides by $F(\psi; t,N;c)$.
Then, using \thetag{\tussenstapje} and analytic continuation to replace
$u$ by $-e^{-i\psi}$, we see that we end up with
$$
\eqalign{
R_{l,m}^{(\nu)}(z) = \sum_{r=0}^l \sum_{s=0}^m &C_{l,m;r,s}^{(\n)} (-1)^{r-s}
t^{\hf (r+s)} c^{r+s}\, R_{l-r,m-s}^{(\nu +r+s)}(\sqrt{1-tc})\cr 
&\times R_{l-r,m-s}^{(\nu +r+s)}(\sqrt{1-c})\,
R_{r,s}^{(\nu -1)}(\sqrt{1-c^{N-1}})\, (-e^{-i\psi})^{s-r}\cr}
$$
where now $z= \sqrt{(1-c)(1-tc)} + c \sqrt{t(1-c^{N-1})}e^{i\psi}$.
Upon substituting $\sqrt{1-c} =\cos\theta_1$, $\sqrt{1-tc} =\cos\theta_2
\, , \, \sqrt{1-c^{N-1}} = \r$ and invoking the homogeneity argument given
in section 7, we obtain the classical identity \thetag{\classadd}.

\bigskip

We end with a few words on the limit case $q\uparrow 1$ of the product formula
\thetag{\prodformlidentity}. As was mentioned in section 7 we apply
\cite{\AsKoo, Thm.~2} to the orthonormal polynomials $p_z(x;n) =
\hat K_z (xc^{-N}; t, nzN; c^{1/nz})$. From \thetag{\Kraworth} we derive that
their orthogonality is given by
$$
\eqalign{
\sum_{x=0}^{nzN} &\frac{(c^N;c^{-1/nz})_{nzN-x} (tc^{1/nz}; c^{1/nz})_{nzN-x}}
{(c^{1/nz}; c^{1/nz})_{nzN-x}}\, (tc^{1/nz})^x \cr
&\times\hat K_z (c^{x/nz}c^{-N}; t, nzN; c^{1/nz})\,
\hat K_{z\pr} (c^{x/nz}c^{-N}; t, nzN; c^{1/nz}) =\d_{zz\pr}.\cr}
$$
So the measure of orthogonality for these polynomials has as support
the set $\{ c^{x/nz} \mid x=0,1,\hdots, nzN\}$ which is contained in
the interval $[0,1]$. Again it is easy to check that for all $k\in\Z$
$$
\eqalign{
&\lim_{n\to\infty} \hat a_{nz+k}(t,nzN+a;c^{1/nz}) = A\cr
&\lim_{n\to\infty} \hat b_{nz+k}(t,nzN+a;c^{1/nz}) = B,\cr}
$$
with $A$ and $B$ as before. This proves that the conditions of 
\cite{\AsKoo, Thm.~2} are all satisfied.

%%%%%%%%%%%%%%%%%References%%%%%%%%%%%%%%%%%%%%%%%%%%%%%%%%%%%%%%%%%%%%%%

\enddocument